\documentclass[11pt,a4paper]{article}
\usepackage{amssymb}

\usepackage{graphicx}
\usepackage{amsmath}
\usepackage{makeidx}
\usepackage{graphicx}
\usepackage{graphicx}
\usepackage{indentfirst}

\newtheorem{result}{Result}[section]

\newcounter{resultnum}[section]
\newtheorem{conclusion}{Conclusion}[section]

\newcounter{conclusionnum}[section]

\newcounter{conditionnum}[section]

\newcounter{conjecturenum}[section]

\newcounter{examplenum}[section]

\newcounter{exercisenum}[section]

\newcounter{lemmanum}[section]

\newcounter{notationnum}[section]
\newtheorem{theorem}{Theorem}[section]

\newcounter{theoremnum}[section]
\newtheorem{definition}{Definition}[section]

\newcounter{definitionnum}[section]
\newtheorem{corollary}{Corollary}[section]

\newcounter{corollarynum}[section]
\newtheorem{remark}{Remark}[section]

\newcounter{remarknum}[section]
\newtheorem{proposition}{Proposition}[section]

\newcounter{propositionnum}[section]

\newcounter{acknowledgementnum}[section]

\newcounter{algorithmnum}[section]

\newcounter{axiomnum}[section]

\newcounter{casenum}[section]

\newcounter{claimnum}[section]

\newcounter{summarynum}[section]

\newcounter{problemnum}[section]
\newenvironment{proof}[1][]{\textbf{Proof.} }{}

\begin{document}

\title{Nonlinear Connections on Gerbes, \\
Clifford--Finsler Modules, and\\
the Index Theorems }
\author{Sergiu I. Vacaru\thanks{%
sergiu$_{-}$vacaru@yahoo.com, svacaru@fields.utoronto.ca } \quad and\ Juan
F. Gonz\'alez--Hern\'andez \thanks{
juanfrancisco.gonzalez@titulado.uam.es, jfgh.teorfizikisto@gmail.com} \\
{\ } \\
\textsl{The Fields Institute for Research in Mathematical Science} \\
\textsl{222 College Street, 2d Floor, } \textsl{Toronto \ M5T 3J1, Canada} \\
and \\
\textsl{Faculty of Mathematics, University "Al. I. Cuza" 
Ia\c si}, \\
\textsl{ 700506, Ia\c si, Romania}}
\date{September 28, 2008}
\maketitle

\begin{abstract}
The geometry of nonholonomic bundle gerbes, provided with nonlinear
connection structure, and nonholonomic gerbe modules is elaborated as the
theory of Clifford modules on nonholonomic manifolds which positively fail
to be spin. We explore an approach to such nonholonomic Dirac operators and
derive the related Atiyah--Singer index formulas. There are considered
certain applications in modern gravity and geometric mechanics of
Clifford--Lagrange/ Finsler gerbes and their realizations as nonholonomic
Clifford and Riemann--Cartan modules.

\vskip0.3cm

\textbf{Keywords:}\ Nonholonomic gerbes, nonlinear connections,
Riem\-ann--Cartan and Lagrange--Finsler spaces, nonholonomic spin structure,
Clifford modules, Dirac operators, the Atiyah--Singer index formulas.

2000 AMS Subject Classification:\

55R65, 53C05, 53C27, 57R20, 57R22, 53B20, 70H99, 81T13, 83C60
\end{abstract}



\section{ Introduction}

In this paper, we elaborate an approach to the Atiyah--Singer index formulas
for manifolds and bundle spaces provided with nonlinear connection (in
brief, N--connection) structure. We follow the methods developed for bundle
gerbes \cite{murray} and bundle gerbe modules \cite{bcnns} and related
results from \cite{aris,mursin}. It should be noted that bundles and gerbes
and their higher generalizations ($n$--gerbes) can be described both in two
equivalent forms:\ in local geometry, with local functions and forms, and in
non--local geometry, by using holonomies and parallel transports, \cite%
{bar91,cp94,mp}, see review \cite{brm}.

The nonholonomic structure of a so--called N--anholonomic space (see, for
instance, \cite{vclalg} and \cite{esv}) is stated by a non--integrable, i.e.
nonholonomic (equivalently, anholonomic), distribution defining a
N--connection structure\footnote{%
the rigorous definitions and notations are given below; for our purposes, it
will be enough to consider nonholonomic spaces defined by N--connections
structures (in brief, called N--anholonomic manifolds)}. Nonholonomic
geometric configurations are naturally derived in modern gravity and string
theories by using generic off--diagonal metrics, generalized connections and
nonholonomic frame structures \cite{vhep2,vs}. The approach can be
elaborated in general form by unifying the concepts of Riemann--Cartan and
Finsler--Lagrange spaces \cite{vggmaf} and their generalizations on gerbes.
It is also related to modelling gravitational field interactions and the
Lagrange and/or Hamilton mechanics and further developments to quantum
deformations \cite{kon,esv}, noncommutative geometry and gravity with
N--anholonomic structures \cite{vncl,vncgg,vesnc}, supermanifolds provided
with N--connection structure \cite{vncsup,bej,az}, as well to nonholonomic
Lie and Clifford algebroids and their applications in constructing new
classes of exact solutions \cite{vesalg}.

The general nonholonomic manifolds fail to be spin and there are substantial
difficulties in definition of curvature which can be revised and solved in
the theory of nonholonomic gerbes. Some of such constructions are relevant
to anomalies in quantum field theory \cite{cmm} when the obstruction to
existence of spin structure is regarded as an anomaly in the global
definition of spinor fields. The typical solution of this problem is to
introduce some additional fields which also have an anomaly in their global
definition but choose a such configuration when both anomalies cancel each
another.

The failure in existence of usual spin structure is differently treated for
the nonholonomic manifolds. At least for the N--anholonomic spaces, it is
possible to define the curvature, which is not a trivial construction for
general nonholonomic manifolds, and the so--called N--adapted (nonholonomic)
Clifford structures with nontrivial N--connection. This problem was firstly
solved for the Finsler--Lagrange spaces \cite{vfs} and their higher order
generalizations \cite{vhs} but it can be also generalized to noncommutative
geometry and gravity models with nontrivial nonholonomic structures and
Lie--Clifford algebroid symmetries, see reviews and recent results in Refs. %
\cite{vstav,vv,vncl}.

This work is devoted to the index theorems for nonholonomic gerbes and
bundle gerbe modules adapted to the N--connection structure. The key idea is
to consider the so--called N--anholonomic spin gerbe defined for any
N--anholonomic manifold (such a nonholonomic gerbe is a usual spin gerbe for
a vanishing N--connection structures and becomes trivial if the basic
manifold is spin). We shall construct ''twisted'' Dirac d--operators
\footnote{%
In brief, we shall write ''d--operators and d--objects'' for operators and
objects distinguished by a N--connection structure, see next section.} and
investigate their properties for Clifford N--anholonomic modules. Then, we
shall define the Chern character of such of such modules and show that the
usual index formula holds for such a definition but being related to the
N--anholonomic structure. The final aim, the proof of index theorems for
various types of N--anholonomic spaces, follows from matching up the
geometric formalism of Clifford modules and nonhlonomic frames with
associated N--connections.

The structure of the paper is as follows:

Section 2 contains an introduction to the geometry of N--anholonom\-ic
manifolds. There are given two equivalent definitions of N--connections,
considered basic geometric objects characterizing them and defined and
computed in abstract form the torsions and curvatures of N--anholonomic
manifolds.

Section 3 is devoted to a study of two explicit examples of N--anholonomic
manifolds: the Lagrange--Finsler spaces and Riemann--Cartan manifolds
provided with N--connection structures. There are proved two main results:
Result \ref{result1}: any regular Lagrange mechanics theory, or Finsler
geometry, can be canonically modelled as a N--anholonomic Riemann--Cartan
manifold with the basic geometric structures (the N--connection, metric and
linear connection) being defined by the fundamental Lagrange, or Finsler,
function; Result \ref{result2}: There are N--anholonomic Einstein--Cartan
(in particular cases, Einstein) spaces parametrized by nontrivial
N--connection structure, nonholonomic frames and, in general, non--Riemann
connections defined as generic off--diagonal solutions in modern gravity.

In section 4, there are considered the lifting of N--anholonomic bundle
gerbes and definition of connection and curvatures on such spaces. We define
the (twisted--anholonomic) Chern characters for bundle gerbes and modules
induced by N--connections and distinguished metric and linear connection
strucutres. We conclude with two important results/ applications of the
theory of nonholonomic gerbes: Result \ref{result3}: Any regular Lagrange
(Finsler) configuration is topologically characterized by its Chern
character computed by using canonical connections defined by the Lagrangian
(fundamental Finsler function). Result \ref{result4}: the geometric
constructions for a N--anholonomic Riemann--Cartan manifold (including exact
solutions in gravity) can be globalized to N--anholonomic gerbe
configurations and characterized by the corresponding Chern character.

Section 5 presents the main result of this paper: The twisted index formula
for N--anholonomic Dirac operators and related gerbe constructions are
stated by Theorem \ref{thmr}). We introduce Clifford d--algebras on
N--anholonomic bundles and define twisted nonholonomic Dirac operators on
Clifford gerbes. We conclude that there are certain fundamental topological
characteristics derived from a regular fundamental Lagrange, or Finsler,
function and that such indices classify new classes of exact solutions in
gravity globalized on gravitational gerbe configurations.

The Appendix contains a set of component formulas for N--connections, metric
and linear connection structures and related torsions and curvatures on
N--anholonomic manifolds. They may considered for some local proofs of
results in the main part of the paper, as well for some applications in
modern physics.

\section{N--Anholonomic Manifolds}

We formulate a coordinate free introduction into the geometry of
nonholonomic manifolds. The reader may consult details in Refs. \cite%
{vhep2,vs,vggmaf,esv}. Here we note that there is a comprehensive study of
nonholonomic and (for integrable structure) of fibred structures in Ref. %
\cite{bejfar1,bejfar2} following the so--called Schouten -- Van Kampen \cite%
{sch} and Vr\v{a}nceanu connections \cite%
{vranceanu1,vranceanu2,vranceanu3,vranceanu4}. Different directions in the
geometry of nonholonomic manifolds were developed for different geometric
structures \cite{mnh1,mnh2,mnh3}, in the geometry of Finsler and Lagrange
spaces \cite{mnh4,mnh5,vaismnh} with applications to mechanics and modern
geometry \cite{versh,leites}. Even from formal point of view all geometric
structures on nonholonomic bundle spaces were rigorously investigated by the
R. Miron's school in Romania, various purposes and applications in modern
physics requested a different class of nonholonomic manifolds with
supersymmetric, noncommutative, Lie algebroid, gerbe etc generalizations %
\cite{vrevfg}. In our approaches, we use such linear and nonlinear
connection structure which can be derived naturally as exact solutions in
modern gravity theories and from certain Lagrangians/ Hamiltonians in the
case of geometric mechanics. Some important component/coordinate formulas
are given in the Appendix.

\subsection{Nonlinear connection structures}

Let $\mathbf{V}$ be a smooth manifold of dimension $(n+m)$ with a local
fibred structure. Two important particular cases are those of a vector
bundle, when we shall write $\mathbf{V=E}$ (with $\mathbf{E}$ being the
total space of a vector bundle $\pi :$ $\mathbf{E}\rightarrow M$ with the
base space $M)$ and of a tangent bundle when we shall consider $\mathbf{V=TM.%
}$\ The differential of a map $\pi :\mathbf{V}\rightarrow M$ defined by
fiber preserving morphisms of the tangent bundles $T\mathbf{V}$ and $TM$ is
denoted by $\pi ^{\top }:T\mathbf{V}\rightarrow TM.$ The kernel of $\pi
^{\top }$ defines the vertical subspace $v\mathbf{V}$ with a related
inclusion mapping $i:v\mathbf{V}\rightarrow T\mathbf{V}.$

\begin{definition}
\label{dnc}A nonlinear connection (N--connection) $\mathbf{N}$ on a manifold
$\mathbf{V}$ is defined by the splitting on the left of an exact sequence%
\begin{equation}
0\rightarrow v\mathbf{V}\overset{i}{\rightarrow }T\mathbf{V}\rightarrow T%
\mathbf{V}/v\mathbf{V}\rightarrow 0,  \label{exseq}
\end{equation}%
i. e. by a morphism of submanifolds $\mathbf{N:\ \ }T\mathbf{V}\rightarrow v%
\mathbf{V}$ such that $\mathbf{N\circ i}$ is the unity in $v\mathbf{V}.$
\end{definition}

The exact sequence (\ref{exseq}) states a nonintegrable (nonholonomic,
equivalently, anholonomic) distribution on $\mathbf{V},$ i.e. this manifold
is nonholonomic. We can say that a N--connection is defined by a global
splitting into conventional horizontal (h) subspace, $\left( h\mathbf{V}%
\right) ,$ and vertical (v) subspace, $\left( v\mathbf{V}\right) ,$
corresponding to the Whitney sum
\begin{equation}
T\mathbf{V}=h\mathbf{V}\oplus _{N}v\mathbf{V}  \label{whitney}
\end{equation}%
where $h\mathbf{V}$ is isomorphic to $M.$ We put the label $N$ to the symbol
$\oplus $ in order to emphasize that such a splitting is associated to a
N--connection structure. In this paper, we shall omit local coordinate
considerations.

For convenience, in Appendix, we give some important local formulas (see,
for instance, the local representation for a N--connection (\ref{nclf})) for
the basic geometric objects and formulas on spaces provided with
N--connection structure. Here, we note that the concept of N--connection
came from E. Cartan's works on Finsler geometry \cite{cart} (see a detailed
historical study in Refs. \cite{ma,esv,vncl} and alternative approaches
developed by using the Ehressmann connection \cite{ehr,dl}). Any manifold
admitting an exact sequence of type (\ref{exseq}) admits a N--connection
structure. If $\mathbf{V=E,}$ a N--connection exists for any vector bundle $%
\mathbf{E}$ over a paracompact manifold $M,$ see proof in Ref. \cite{ma}.

The geometric objects on spaces provided with N--connection structure are
denoted by ''bolfaced'' symbols. Such objects may be defined in
''N--adapted'' form by considering h-- and v--decompositions (\ref{whitney}%
). Following conventions from \cite{ma,vfs,vstav,vncl}, one call such
objects to be d--objects (i. e. they are distinguished by the N--connection;
one considers d--vectors, d--forms, d--tensors, d--spinors, d--connections,
....). For instance, a d--vector is an element $\mathbf{X}$ of the module of
the vector fields $\chi (\mathbf{V)}$ on $\mathbf{V},$ which in N--adapted
form may be written
\begin{equation*}
\mathbf{X=}\ h\mathbf{X}+v\mathbf{X}\mbox{ or }\mathbf{X=\ }X\oplus _{N}\
^{\bullet }X,
\end{equation*}%
where $h\mathbf{X}$ (equivalently, $X$) is the h--component and $v\mathbf{X}$
(equivalently, $\ ^{\bullet }X)$ is the v--component of $\mathbf{X.}$

A N--connection is characterized by its \textbf{N--connection curvature }%
(the Nijenhuis tensor)%
\begin{equation}
\Omega (\mathbf{X,Y})\doteqdot \lbrack \ ^{\bullet }X,\ ^{\bullet }Y]+\
^{\bullet }[\mathbf{X,Y}]-\ ^{\bullet }[\ ^{\bullet }X\mathbf{,Y}]-\
^{\bullet }[\mathbf{X,}\ ^{\bullet }Y]  \label{njht}
\end{equation}%
for any $\mathbf{X,Y\in }\chi (\mathbf{V),}$ where $[\mathbf{X,Y]\doteqdot
XY-}$ $\mathbf{YX}$ and $\ ^{\bullet }[\mathbf{,]}$ is the v--projection of $%
[\mathbf{,],}$ see also the coordinate formula (\ref{ncurv}) in Appendix.
This d--object $\Omega $ was introduced in Ref. \cite{grifone} in order to
define the curvature of a nonlinear connection in the tangent bundle over a
smooth manifold. But this can be extended for any nonholonomic manifold,
nonholnomic Clifford structure and any noncommutative / supersymmetric
versions of bundle spaces provided with N--connection structure, i. e. with
nonintegrable distributions of type (\ref{whitney}), see \cite{esv,vncl,vv}.

\begin{proposition}
\label{pddv}A N--connection structure on $\mathbf{V}$ defines a nonholonomic
N--adapted frame (vielbein) structure $\mathbf{e}=(e,^{\bullet }e)$ and its
dual $\widetilde{\mathbf{e}}=\left( \widetilde{e},\ ^{\bullet }\widetilde{e}%
\right) $ with $e$ and $\ ^{\bullet }\widetilde{e}$ linearly depending on
N--connection coefficients.\
\end{proposition}

\begin{proof}
It follows from explicit local constructions, see formulas (\ref{ddif}), (%
\ref{dder}) and (\ref{anhrel})\ in Appendix.$\square $
\end{proof}

\begin{definition}
A manifold \ $\mathbf{V}$ is called N--anholonomic if it is defined a local
(in general, nonintegrable) distribution (\ref{whitney}) on its tangent
space $T\mathbf{V,}$ i.e. $\mathbf{V}$ is N--anholonomic if it is enabled
with a N--connection structure (\ref{exseq}).
\end{definition}

All spinor and gerbe constructions in this paper will be performed for
N--anholonomic manifolds.

\subsection{Curvatures and torsions of N--anholonomic manifolds}

One can be defined N--adapted linear connection and metric structures on $%
\mathbf{V:}$

\begin{definition}
\label{ddc}A distinguished connection (d--connection) $\mathbf{D}$ on a
N--anho\-lo\-no\-mic manifold $\mathbf{V}$ is a linear connection conserving
under parallelism the Whitney sum (\ref{whitney}). For any $\mathbf{X\in }%
\chi (\mathbf{V),}$ one have a decomposition into h-- and v--covariant
derivatives,%
\begin{equation}
\mathbf{D}_{\mathbf{X}}\mathbf{\doteqdot X}\rfloor \mathbf{D=\ }X\rfloor
\mathbf{D+}\ ^{\bullet }X\rfloor \mathbf{D=}D_{X}+\ ^{\bullet }D_{X}.
\label{dconcov}
\end{equation}
\end{definition}

The symbol ''$\rfloor "$ in (\ref{dconcov}) denotes the interior product. We
shall write conventionally that $\mathbf{D=}(D,\ ^{\bullet }D).$

For any d--connection $\mathbf{D}$ on a N--anholonomic manifold $\mathbf{V},$
it is possible to define the curvature and torsion tensor in usual form but
adapted to the Whitney sum (\ref{whitney}):

\begin{definition}
The torsion
\begin{equation}
\mathbf{T(X,Y)\doteqdot D}_{\mathbf{X}}\mathbf{Y-D}_{\mathbf{Y}}\mathbf{%
X-[X,Y]}  \label{tors1}
\end{equation}%
of a d--connection $\mathbf{D=}(D,\ ^{\bullet }D)\mathbf{,}$ for any $%
\mathbf{X,Y\in }\chi (\mathbf{V),}$ has a N--adapted decomposition
\begin{equation}
\mathbf{T(X,Y)=T(}X,Y\mathbf{)+T(}X,\ ^{\bullet }Y\mathbf{)+T(\ }^{\bullet
}X,Y\mathbf{)+T(\ }^{\bullet }X,\ ^{\bullet }Y\mathbf{).}  \label{tors2}
\end{equation}
\end{definition}

By further h- and v--projections of (\ref{tors2}), denoting $h\mathbf{%
T\doteqdot }T$ and $v\mathbf{T\doteqdot \ }^{\bullet }T,$ taking in the
account that $h\mathbf{[\ }^{\bullet }X,\ ^{\bullet }Y\mathbf{]=}0\mathbf{,}$
one proves

\begin{theorem}
\label{tht}The torsion of a d--connection $\mathbf{D=}(D,\ ^{\bullet }D)$ is
defined by five nontrivial d--torsion fields adapted to the h-- and
v--splitting by the N--connection structure%
\begin{eqnarray*}
T(X,Y) & \doteqdot &D_{X}Y-D_{Y}X-h[X,Y], \\
\mathbf{\ }^{\bullet } T(X,Y) &\doteqdot & {\ }^{\bullet }[Y,X], \\
T(X,\mathbf{\ }^{\bullet }Y) &\mathbf{\doteqdot }&\mathbf{-\ }^{\bullet
}D_{Y}X-h[X,\mathbf{\ }^{\bullet }Y], \\
\mathbf{\ }^{\bullet }T(X,\mathbf{\ }^{\bullet }Y) &\doteqdot& {\ }^{\bullet
}D_{X}Y-\mathbf{\ }^{\bullet }[X,\mathbf{\ }^{\bullet }Y], \\
\mathbf{\ }^{\bullet }T(\mathbf{\ }^{\bullet }X,\mathbf{\ }^{\bullet }Y) &%
\mathbf{\doteqdot }&\mathbf{\ }^{\bullet }D_{X}\mathbf{\ }^{\bullet }Y-%
\mathbf{\ }^{\bullet }D_{Y}\mathbf{\ }^{\bullet }X-\mathbf{\ }^{\bullet }[%
\mathbf{\ }^{\bullet }X,\mathbf{\ }^{\bullet }Y].
\end{eqnarray*}%
\
\end{theorem}

The d--torsions $T(X,Y),\mathbf{\ }^{\bullet }T(\mathbf{\ }^{\bullet }X,%
\mathbf{\ }^{\bullet }Y)$ are called respectively the $h(hh)$--torsion, $%
v(vv)$--torsion and so on. The formulas (\ref{dtors}) in Appendix present a
local proof of this Theorem.

\begin{definition}
The curvature of a d--connection $\mathbf{D=}(D,\ ^{\bullet }D)$ is defined%
\textbf{\ }
\begin{equation}
\mathbf{R(X,Y)\doteqdot D}_{\mathbf{X}}\mathbf{D}_{\mathbf{Y}}\mathbf{-D}_{%
\mathbf{Y}}\mathbf{D}_{\mathbf{X}}\mathbf{-D}_{\mathbf{[X,Y]}}  \label{curv1}
\end{equation}
for any $\mathbf{X,Y\in }\chi (\mathbf{V).}$
\end{definition}

Denoting $h\mathbf{R}=R$ and $v\mathbf{R}=\ ^{\bullet }R,$ by
straightforward calculations, one check the properties%
\begin{eqnarray*}
R(\mathbf{X,Y})\mathbf{\ }^{\bullet }Z &=&0,\mathbf{\ }^{\bullet }R(\mathbf{%
X,Y})Z=0, \\
\mathbf{R(X,Y)Z} &\mathbf{=}&R\mathbf{(X,Y)}Z+\mathbf{\ }^{\bullet }R\mathbf{%
(X,Y)\ }^{\bullet }Z
\end{eqnarray*}%
for any for any $\mathbf{X,Y,Z\in }\chi (\mathbf{V).}$

\begin{theorem}
\label{thr}The curvature $\mathbf{R}$ of a d--connection $\mathbf{D=}(D,\
^{\bullet }D)$ is completely defined by six d--curvatures
\begin{eqnarray*}
\mathbf{R(}X\mathbf{,}Y\mathbf{)}Z &\mathbf{=}&\left(
D_{X}D_{Y}-D_{Y}D_{X}-D_{[X,Y]}-\mathbf{\ }^{\bullet }D_{[X,Y]}\right) Z, \\
\mathbf{R(}X\mathbf{,}Y\mathbf{)\ }^{\bullet }Z &\mathbf{=}&\left(
D_{X}D_{Y}-D_{Y}D_{X}-D_{[X,Y]}-\mathbf{\ }^{\bullet }D_{[X,Y]}\right) \
^{\bullet }Z, \\
\mathbf{R(\ }^{\bullet }X\mathbf{,}Y\mathbf{)}Z &\mathbf{=}&\left( \mathbf{\
}^{\bullet }D_{X}D_{Y}-D_{Y}\mathbf{\ }^{\bullet }D_{X}-D_{[\mathbf{\ }%
^{\bullet }X,Y]}-{}^{\bullet }D_{[{}^{\bullet}X,Y]}\right)Z, \\
\mathbf{R}({\ }^{\bullet}X,Y){\ }^{\bullet }Z &=&\left({\ }^{\bullet} D_{X}{%
\ }^{\bullet }D_{Y}-{\ }^{\bullet}D_{Y}{\ }^{\bullet}D_{X}- D_{[{\ }%
^{\bullet}X,Y]}-{}^{\bullet}D_{[{\ }^{\bullet}X,Y]}\right)^{\bullet}Z, \\
\mathbf{R(\ }^{\bullet }X\mathbf{,\ }^{\bullet }Y\mathbf{)}Z &=&\left(
\mathbf{\ }^{\bullet }D_{X}D_{Y}-D_{Y}\mathbf{\ }^{\bullet }D_{X}-\mathbf{\ }%
^{\bullet }D_{[\mathbf{\ }^{\bullet }X,\mathbf{\ }^{\bullet }Y]}\right) Z, \\
\mathbf{R}(^{\bullet }X{,}^{\bullet }Y){}^{\bullet}Z &=& \left( \mathbf{\ }%
^{\bullet }D_{X}D_{Y}-D_{Y}\mathbf{\ }^{\bullet}D_{X}-\mathbf{\ }%
^{\bullet}D_{[{\ }^{\bullet }X,{\ }^{\bullet }Y]}\right){\ }^{\bullet}Z.
\end{eqnarray*}
\end{theorem}

The proof of Theorems \ref{tht} and \ref{thr} is given for vector bundles
provided with N--connection structure in Ref. \cite{ma}. Similar Theorems
and respective proofs hold true for superbundles \cite{vncsup}, for
noncommutative projective modules \cite{vncl} and for N--anholonomic
metric--affine spaces \cite{vggmaf}, where there are also give the main
formulas in abstract coordinate form. The formulas (\ref{dcurv}) from
Appendix consist a coordinate proof of Theorem \ref{thr}.

\begin{definition}
A metric structure $\ \breve{g}$ on a N--anholonomic space $\mathbf{V}$ is a
symmetric covariant second rank tensor field which is not degenerated and of
constant signature in any point $\mathbf{u\in V.}$
\end{definition}

In general, a metric structure is not adapted to a N--connection structure.

\begin{definition}
\label{ddm}A d--metric $\mathbf{g}=g\oplus _{N}\ ^{\bullet }g$ is a usual
metric tensor which contracted to a d--vector results in a dual d--vector,
d--covector (the duality being defined by the inverse of this metric tensor).
\end{definition}

The relation between arbitrary metric structures and d--metrics is
established by

\begin{theorem}
\label{tdm}Any metric $\ \breve{g}$ can be equivalently transformed into a
d--metric
\begin{equation}
\mathbf{g}=g(X,Y)+\mathbf{\ }^{\bullet }g(\mathbf{\ }^{\bullet }X,\mathbf{\ }%
^{\bullet }Y)  \label{dmetra}
\end{equation}%
for a corresponding N--connection structure.
\end{theorem}

\begin{proof}
We introduce denotations $h\breve{g}(X,Y)\doteqdot g(X,Y)$ and $\mathbf{\ }v%
\breve{g}(\mathbf{\ }^{\bullet }X,\mathbf{\ }^{\bullet }Y)$ $=\mathbf{\ }%
^{\bullet }g(\mathbf{\ }^{\bullet }X,\mathbf{\ }^{\bullet }Y)$ and try to
find a N--connection when
\begin{equation}
\breve{g}(X,\mathbf{\ }^{\bullet }Y)=0  \label{algn01}
\end{equation}%
for any $\mathbf{X,Y\in }\chi (\mathbf{V).}$ In local form, \ the equation (%
\ref{algn01}) is just an algebraic equation for $\mathbf{N}=\{N_{i}^{a}\},$
see formulas (\ref{metr}), (\ref{ansatz}) and (\ref{dmetr}) and related
explanations in Appendix. $\square $
\end{proof}

\begin{definition}
A d--connection $\mathbf{D}$ on $\mathbf{V}$ is said to be metric, i.e. it
satisfies the metric compatibility (equivalently, metricity) conditions with
a metric $\ \breve{g}$ and its equivalent d--metric $\mathbf{g},$ if there
are satisfied the conditions
\begin{equation}
\mathbf{D}_{\mathbf{X}}\mathbf{g=0.}  \label{mcc}
\end{equation}
\end{definition}

Considering explicit h-- and v--projecting of (\ref{mcc}), one proves

\begin{proposition}
\label{pmcc}A d--connection $\mathbf{D}$ on $\mathbf{V}$ is metric if and
only if
\begin{equation*}
D_{X}g=0,\ D_{X}\mathbf{\ }^{\bullet }g=0,\ ^{\bullet }D_{X}g=0,\ ^{\bullet
}D_{X}\mathbf{\ }^{\bullet }g=0.
\end{equation*}
\end{proposition}

One holds this important

\begin{conclusion}
Following Propositions \ref{pddv} and \ref{pmcc} and Theorem \ref{tdm}, we
can elaborate the geometric constructions on a N--anholonomic manifold $%
\mathbf{V}$ $\ $\ in N--adapted form by considering N--adapted frames $%
\mathbf{e}=(e,^{\bullet }e)$ and co--frames $\widetilde{\mathbf{e}}=\left(
\widetilde{e},\ ^{\bullet }\widetilde{e}\right) ,$ d--connection $\mathbf{D}$
and d--metric fields $\mathbf{g}=[g,\mathbf{\ }^{\bullet }g].$
\end{conclusion}

In Riemannian geometry, there is a preferred linear Levi--Civita connection $%
\nabla $ which is metric compatible and torsionless, i.e.
\begin{equation*}
\ ^{\nabla }\mathbf{T(X,Y)\doteqdot }\nabla _{\mathbf{X}}\mathbf{Y-}\nabla _{%
\mathbf{Y}}\mathbf{X-[X,Y]=}0,
\end{equation*}%
and defined by the metric structure. On a general N--anholonomic manifold $%
\mathbf{V}$ provided with a d--metric structure $\mathbf{g}=[g,\mathbf{\ }%
^{\bullet }g],$ the Levi--Civita connection defined by this metric is not
adapted to the N--connection, i. e. to the splitting (\ref{whitney}). The
h-- and v--distributions are nonintegrable ones and any d--connection
adapted to a such splitting contains nontrivial d--torsion coefficients.
Nevertheless, one exists a minimal extension of the Levi--Civita connection
to a canonical d--connection which is defined only by a metric $\breve{g}.$

\begin{theorem}
\label{thcdc}For any d--metric $\mathbf{g}=[g,\mathbf{\ }^{\bullet }g]$ on a
N--anholonomic manifold $\mathbf{V,}$ there is a unique metric canonical
d--connection $\widehat{\mathbf{D}}$ satisfying the conditions $\widehat{%
\mathbf{D}}\mathbf{g=}0$ and with vanishing $h(hh)$--torsion, $v(vv)$%
--torsion, i. e. $\widehat{T}(X,Y)=0$ and $\mathbf{\ }^{\bullet }\widehat{T}(%
\mathbf{\ }^{\bullet }X,\mathbf{\ }^{\bullet }Y)=0.$
\end{theorem}

\begin{proof}
The formulas (\ref{cdc}) and (\ref{candcon}) and related discussions in
Appendix give a proof, in component form, of this Theorem.$\square $
\end{proof}

The following Corollary gathers some basic information about
N--anho\-lo\-no\-mic manifolds.

\begin{corollary}
\label{cncs}A N--connection structure defines three important geometric
objects:

\begin{enumerate}
\item a (pseudo) Euclidean N--metric structure $\ ^{\eta }\mathbf{g}=\eta
\oplus _{N}\ ^{\bullet }\eta ,$ i.e. a d--metric with (pseudo) Euclidean
metric coefficients with respect to $\widetilde{\mathbf{e}}$ defined only by
$\mathbf{N;}$

\item a N--metric canonical d--connection $\widehat{\mathbf{D}}^{N}$ defined
only by $^{\eta }\mathbf{g}$ and $\mathbf{N;}$

\item a nonmetric Berwald type linear connection $\mathbf{D}^{B}.$
\end{enumerate}
\end{corollary}

\begin{proof}
Fixing a signature for the metric, $sign^{\ \eta }\mathbf{g=(\pm ,\pm
,...,\pm ),}$ we introduce these values in (\ref{dmetr}) we get $^{\eta }%
\mathbf{g}=\eta \oplus _{N}\ ^{\bullet }\eta $ of type (\ref{dmetra}), i.e.
we prove the point 1. The point 2 is to be proved by an explicit
construction by considering the coefficients of $^{\eta }\mathbf{g}$ into (%
\ref{candcon}). This way, we get a canonical d--connection induced by the
N--connection coefficients and satisfying the metricity conditions (\ref{mcc}%
). In an approach to Finsler geometry \cite{cbs}, one emphasizes the
constructions derived for the so--called Berwald type d--connection $\mathbf{%
D}^{B},$ considered to be the ''most'' minimal (linear on $\Omega )$
extension of the Levi--Civita connection, see formulas (\ref{dbc}). Such
d--connections can be defined for an arbitrary d--metric $\mathbf{g}=[g,%
\mathbf{\ }^{\bullet }g],$ or for any $^{\eta }\mathbf{g}=\eta \oplus _{N}\
^{\bullet }\eta .$ They are only ''partially'' metric because, for instance,
$D^{B}g=0$ and $\mathbf{\ }^{\bullet }D^{B}\mathbf{\ }^{\bullet }g=0$ but,
in general, $D^{B}\mathbf{\ }^{\bullet }g\neq 0$ and $\mathbf{\ }^{\bullet
}D^{B}g\neq 0,$ i. e. $\mathbf{D}^{B}\mathbf{g}\neq 0,$ see Proposition \ref%
{pmcc}. It is a more sophisticate problem to define spinors and
supersymmetric physically valued models for such Finsler spaces, see
discussions in \cite{vncl,vclalg,vggmaf}. $\square $
\end{proof}

\begin{remark}
The geometrical objects $\widehat{\mathbf{D}}^{N},\mathbf{D}^{B}$ for $%
^{\eta }\mathbf{g,}$ nonholonomic bases $\mathbf{e}=(e,^{\bullet }e)$ and $%
\widetilde{\mathbf{e}}=\left( \widetilde{e},\ ^{\bullet }\widetilde{e}%
\right) ,$ see Proposition \ref{pddv} and the N--connection curvature $%
\mathbf{\Omega }$ (\ref{njht}), define completely the main properties of a
N--anholonomic manifold $\mathbf{V.}$
\end{remark}

It is possible to extend the constructions for any additional d--metric and
canonical d--connection structures. For our considerations on nonholnomic
Clifford/spinor structures, the class of metric d--connections plays a
preferred role. That why we emphasize the physical importance of
d--connections $\widehat{\mathbf{D}}$ and $\widehat{\mathbf{D}}^{N}$ instead
of $\mathbf{D}^{B}$ or any other nonmetric d--connections.

Finally, in this section, we note that the d--torsions and d--curvatures on
N--anholonomic manifolds can be computed for any type of d--connection
structure, see Theorems \ref{tht} and \ref{thr} and the component formulas (%
\ref{dtors}) and (\ref{dcurv}).

\section{Examples of N--anholonomic spaces:}

For corresponding parametrizations of the N--connection, d--metric and
d--connec\-ti\-on coefficients of a N--anholonomic space, it is possible to
model various classes of (generalized) Lagrange, Finsler and Riemann--Cartan
spa\-ces. We briefly analyze three such nonholonomic geometric structures.

\subsection{Lagrange--Finsler geometry}

This class of geometries is usually defined on tangent bundles \cite{ma} but
it is possible to model such structures on general N--anholonomic manifolds,
in particular in (pseudo) Riemannian and Riemann--Cartan geometry if
nonholonomic frames are introduced into consideration \cite%
{vhep2,vs,vggmaf,vesalg}. Let us outline the first approach when the
N--anholonomic manifold $\mathbf{V}$ is taken to be just a tangent bundle $%
(TM,\!\pi ,\!M),$ where $M$ is a $n$--dimensional base manifold, $\pi $ is a
surjective projection and $TM$ is the total space. One denotes by $%
\widetilde{TM}=TM\backslash \{0\}$ where $\{0\}$ means the null section of
map $\pi .$

We consider a differentiable fundamental Lagrange function $L(x,y)$ defined
by a map $L:(x,y)\in TM\rightarrow L(x,y)\in \mathbb{R}$ of class $\mathcal{C%
}^{\infty }$ on $\widetilde{TM}$ and continuous on the null section $%
0:M\rightarrow TM$ of $\pi .$ The values $x=\{x^{i}\}$ are local coordinates
on $M$ and $(x,y)=(x^{i},y^{k})$ are local coordinates on $TM.$ For
simplicity, we consider this Lagrangian to be regular, i.e. with
nondegenerated Hessian
\begin{equation}
\ ^{L}g_{ij}(x,y)=\frac{1}{2}\frac{\partial ^{2}L(x,y)}{\partial
y^{i}\partial y^{j}}  \label{lqf}
\end{equation}%
when $rank\left| g_{ij}\right| =n$ on $\widetilde{TM}$ and the left up ''L''
is an abstract label pointing that the values are defined by the Lagrangian $%
L.$

\begin{definition}
\label{dlg}A Lagrange space is a pair $L^{n}=\left[ M,L(x,y)\right] $ with
the tensor $\ ^{L}g_{ij}(x,y)$ being of constant signature over $\widetilde{%
TM}.$
\end{definition}

The notion of Lagrange space was introduced by J. Kern \cite{kern} and
elaborated in details in Ref. \cite{ma} as a natural extension of Finsler
geometry.

\begin{theorem}
\label{tcslg}There are canonical N--connection $\ ^{L}\mathbf{N,}$ almost
complex $^{L}\mathbf{F,}$ d--metric $\ ^{L}\mathbf{g}$ and d--connection $\
^{L}\widehat{\mathbf{D}}$ structures defined by a regular Lagrangian $L(x,y)$
and its Hessian $\ ^{L}g_{ij}(x,y)$ (\ref{lqf}).
\end{theorem}

\begin{proof}
The canonical $\ ^{L}\mathbf{N}$ is defined by certain nonlinear spray
configurations related to the solutions of Euler--Lagrange equations, see
the local formula (\ref{cncl}) in Appendix. It is given there the explicit
matrix representation of $^{L}\mathbf{F}$ (\ref{acs1}) which is a usual
definition of almost complex structure, after $\ ^{L}\mathbf{N}$ and
N--adapted bases have been constructed. The d--metric (\ref{slm}) is a local
formula for $\ ^{L}\mathbf{g.}$ Finally, the canonical d--connection $\ ^{L}%
\widehat{\mathbf{D}}$ is a usual one but for $\ ^{L}\mathbf{g}$ and $\ ^{L}%
\mathbf{N}$ on $\widetilde{TM}.\square $
\end{proof}

A similar Theorem can be formulated and proved for the Finsler geometry:

\begin{remark}
\label{rcsfg}A Finsler space defined by a fundamental Finsler function $%
F(x,y),$ being homogeneous of type $F(x,\lambda y)=\lambda F(x,y),$ for
nonzero $\lambda \in \mathbb{R},$ may be considered as a particular case of
Lagrange geometry when $L=F^{2}.$
\end{remark}

From the Theorem \ref{tcslg} \ and Remark \ref{rcsfg} one follows:

\begin{result}
\label{result1}Any Lagrange mechanics with regular Lagrangian $L(x,y)$ (any
Finsler geometry with fundamental function $F(x,y))$ can be modelled as a
nonhlonomic Riemann--Cartan geometry with canonical structures $\ ^{L}%
\mathbf{N,}$ $\ ^{L}\mathbf{g}$ and $\ ^{L}\widehat{\mathbf{D}}$ ($\ ^{F}%
\mathbf{N,}$ $\ ^{F}\mathbf{g}$ and $\ ^{F}\widehat{\mathbf{D}})$ defined on
a corresponding N--anholonomic manifold $\mathbf{V}.$
\end{result}

It was concluded that any regular Lagrange mechanics/Finsler geometry can be
geometrized/modelled as an almost K\"{a}hler space with canonical
N--connection distribution, see \cite{ma} and, for N--anholonomic Fedosov
manifolds, \cite{esv}. Such approaches based on almost complex structures
are related with standard sympletic geometrizations of classical mechanics
and field theory, for a review of results see Ref. \cite{dl}.

For applications in optics of nonhomogeneous media \cite{ma} and gravity
(see, for instance, Refs. \cite{vhep2,vs,vggmaf,vclalg,vesalg}), one
considers metrics of type $g_{ij}\sim e^{\lambda (x,y)}\ ^{L}g_{ij}(x,y)\ $\
which can not be derived from a mechanical Lagrangian but from an effective
''energy'' function. In the so--called generalized Lagrange geometry, one
introduced Sasaki type metrics (\ref{slm}), see the Appendix, with any
general coefficients both for the metric and N--connection.

\subsection{N--connections and gravity}

Now we show how N--anholonomic configurations can defined in gravity
theories. In this case, it is convenient to work on a general manifold $%
\mathbf{V},\dim \mathbf{V}=n+m$ enabled with a global N--connection
structure, instead of the tangent bundle $\widetilde{TM}.$

For N--connection splittings of (pseudo) Riemannian--Cartan spaces of
dimension $(n+m)$ (there were also considered (pseudo) Riemannian
configurations), the Lagrange and Finsler type geometries were modelled by
N--anholonomic structures as exact solutions of gravitational field
equations \cite{vggmaf,vs,vesnc}. Inversely, all approaches to (super)
string gravity theories deal with nontrivial torsion and (super) vielbein
fields which under corresponding parametrizations model N--anholonomic
spaces \cite{vncsup,vhs,vv}. We summarize here some geometric properties of
gravitational models with nontrivial N--anholonomic structure.

\begin{definition}
A N--anholonomic Riemann--Cartan manifold $\ ^{RC}\mathbf{V}$ is defined by
a d--metric $\mathbf{g}$ and a metric d--connection $\mathbf{D}$ structures
adapted to an exact sequence splitting (\ref{exseq}) defined on this
manifold.
\end{definition}

The d--metric structure $\mathbf{g}$ on$\ ^{RC}\mathbf{V}$ is of type (\ref%
{dmetra}) and satisfies the metricity conditions (\ref{mcc}). With respect
to a local coordinate basis, the metric $\mathbf{g}$ is parametrized by a
generic off--diagonal metric ansatz (\ref{ansatz}), see Appendix. In a
particular case, we can take $\mathbf{D=}\widehat{\mathbf{D}}$ and treat the
torsion $\widehat{\mathbf{T}}$ as a nonholonomic frame effect induced by
nonintegrable N--splitting. For more general applications, we have to
consider additional torsion components, for instance, by the so--called $H$%
--field in string gravity.

Let us denote by $Ric(\mathbf{D})$ and $Sc(\mathbf{D}),$ respectively, the
Ricci tensor and curvature scalar defined by any metric d--connection $%
\mathbf{D}$ and d--metric $\mathbf{g}$ on $\ ^{RC}\mathbf{V,}$ see also the
component formulas (\ref{dricci}), (\ref{sdccurv}) \ and (\ref{enstdt})\ in
Appendix. The Einstein equations are
\begin{equation}
En(\mathbf{D})\doteqdot Ric(\mathbf{D})-\frac{1}{2}\mathbf{g}Sc(\mathbf{D})=%
\mathbf{\Upsilon }  \label{einsta}
\end{equation}%
where the source $\mathbf{\Upsilon }$ reflects any contributions of matter
fields and corrections from, for instance, string/brane theories of gravity.
In a closed physical model, the equation (\ref{einsta}) have to be completed
with equations for the matter fields, torsion contributions and so on (for
instance, in the\ Einstein--Cartan theory one considers algebraic equations
for the torsion and its source)... It should be noted here that because of
nonholonomic structure of $^{RC}\mathbf{V,}$ the tensor $Ric(\mathbf{D})$ is
not symmetric and that $\mathbf{D}\left[ En(\mathbf{D})\right] \neq 0$ which
imposes a more sophisticate form of conservation laws on such spaces with
generic ''local anisotropy'', see discussion in \cite{vggmaf,vstav} (this is
similar with the case when the nonholonomic constraints in Lagrange
mechanics modifies the definition of conservation laws). A very important
class of models can be elaborated when $\mathbf{\Upsilon =}diag\left[
\lambda ^{h}(\mathbf{u})g,\lambda ^{v}(\mathbf{u})\mathbf{\ }^{\bullet }g%
\right] ,$ which defines the so--called N--anholonomic Einstein spaces.

\begin{result}
\label{result2}Various classes of vacuum and nonvacuum exact solutions of (%
\ref{einsta}) parametrized by generic off--diagonal metrics, nonholonomic
vielbeins and Levi--Civita or non--Riemannian connections in Einstein and
extra dimension gravity models define explicit examples of N--anholonomic
Einstein--Cartan (in particular, Einstein) spaces.
\end{result}

Such exact solutions (with noncommutative, algebroid, toro\-id\-al,
ellipsoid, ... symmetries) have been constructed in Refs. \cite%
{vhep2,vs,esv,vncl,vesnc,vesalg,vclalg,vggmaf,vstav}. We note that a
subclass of N--anholonomic Einstein spaces was related to generic
off--diagonal solutions in general relativity by such nonholonomic
constraints when $Ric(\widehat{\mathbf{D}})=Ric(\mathbf{\nabla })$ even $%
\widehat{\mathbf{D}}\neq \nabla ,$ where $\widehat{\mathbf{D}}$ is the
canonical d--connection and $\nabla $ is the Levi--Civita connection, see
formulas (\ref{cdc}) and (\ref{cdca}) in Appendix and details in Ref. \cite%
{vesalg}.

A direction in modern gravity is connected to analogous gravity models when
certain gravitational effects and, for instance, black hole configurations
are modelled by optical and acoustic media, see a recent review or results
in \cite{blv}. Following our approach on geometric unification of gravity
and Lagrange regular mechanics in terms of N--anholonomic spaces, one holds

\begin{theorem}
\label{tlr}A Lagrange (Finsler) space can be canonically modelled as an
exact solution of the Einstein equations (\ref{einsta}) on a N--anholonomic
Riemann--Cartan space if and only if the canonical N--connection $\ ^{L}%
\mathbf{N} $ ($\ ^{F}\mathbf{N}$)$\mathbf{,}$ d--metric $\ ^{L}\mathbf{g}$ ($%
\ ^{F}\mathbf{g)}$ and d--connection $\ ^{L}\widehat{\mathbf{D}}$ ($\ ^{F}%
\widehat{\mathbf{D}})$ $\ $structures defined by the corresponding
fundamental Lagrange function $L(\mathbf{x,y})$ (Finsler function $F(\mathbf{%
x,y}))$ satisfy the gravitational field equations for certain physically
reasonable sources $\mathbf{\Upsilon .}$
\end{theorem}

\begin{proof}
We sketch the idea: It can be performed in local form by considering the
Einstein tensor (\ref{enstdt}) defined by the $\ ^{L}\mathbf{N}$ ($\ ^{F}%
\mathbf{N}$) in the form (\ref{cncl}) and $\ ^{L}\mathbf{g}$ ($\ ^{F}\mathbf{%
g)}$ in the form (\ref{slm}) \ inducing the canonical d--connection $\ ^{L}%
\widehat{\mathbf{D}}$ ($\ ^{F}\widehat{\mathbf{D}}).$ For certain zero or
nonzero $\mathbf{\Upsilon }$, such N--anholonomic configurations may be
defined by exact solutions of the Einstein equations for a d--connection
structure. A number of explicit examples were constructed for N--anholonomic
Einstein spaces \cite{vhep2,vs,esv,vncl,vesnc,vesalg,vclalg,vggmaf,vstav}.$%
\square $
\end{proof}

It should be noted that Theorem \ref{tlr} states explicit conditions when
the Result \ref{result1} holds for N--anholonomic Einstein spaces.

\begin{conclusion}
Generic off--diagonal metric and vielbein structures in gra\-vi\-ty and
regular Lagrange mechanics models can be geometrized in a unified form on
N--anholonomic manifolds. In general, such spaces are not spin and this
presents a strong motivation for elaborating the theory of nonholonomic
gerbes and related Clifford/spinor structures developed in this work.
\end{conclusion}

Following this Conclusion, it is not surprizing that a lot of gravitational
effects (black hole configurations, collapse scenaria, cosmological
anisotropi\-es etc) can be modelled in nonlinear fluid, acoustic or optic
media.

\section{Lifts of Nonholonomic Bundle Gerbes and Connections}

In this section, we present an introduction into the geometry of lifts of
nonholonomic bundle gerbes and related N--anholonomic modules. We define
connections and curvatures for such bundle modules. This material
reproduces, in the corresponding holonomic limits, certain fundamental
results from \cite{murray,bcnns,mursin,aris}.

\subsection{N--anholonomic bundle gerbes and their lifts}

\subsubsection{Local constructions}

On N--anholonomic manifolds, one deals with nonintegrable h-- and
v--splitting of geometric objects, described by the so--called d--objects
(for instance, d--vectors, d--spinors, d--tensors, d--connections, ... like
we considered in the previous section). It is convenient to introduce the
concept of Lie d--group $\mathbf{G}=(G,\ ^{\bullet }G)$ \cite%
{vfs,vhs,vncl,vstav,vv,vclalg} which is just a couple of two usual Lie
groups $G$ and$\ ^{\bullet }G$ associated to a N--connection splitting (\ref%
{whitney}). \ We conventionally consider a central extension of a finite
dimensional of Lie d--groups $\mathbf{G}$ to $\mathbf{\check{G},}$ defined
by a map $\pi :$ $\mathbf{\check{G}\rightarrow G}$ such that it is defined
the exact sequence
\begin{equation}
0\rightarrow \mathbb{Z}_{k}\rightarrow \mathbf{\check{G}\rightarrow
G\rightarrow }1  \label{dextseq}
\end{equation}%
where $\mathbb{Z}_{k}=\mathbb{Z}/k\mathbb{Z}$ denotes the cyclic subgroup of
the circle $U(1).$ This sequence of d--groups splits into respective
horizontal component%
\begin{equation*}
0\rightarrow \mathbb{Z}_{k}\rightarrow \check{G}\mathbf{\rightarrow }G%
\mathbf{\rightarrow }1
\end{equation*}%
and vertical component%
\begin{equation*}
0\rightarrow \mathbb{Z}_{k}\rightarrow \ ^{\bullet }\check{G}\mathbf{%
\rightarrow }\ ^{\bullet }G\mathbf{\rightarrow }1.
\end{equation*}

Let us denote by $\mathbf{U}$ and $\mathbf{\check{U}}$ the corresponding
right principal sets: the are just $\mathbf{G}$ and $\mathbf{\check{G}}$ but
conventionally re--defined in order to consider distinguished (not mixing
the h-- and v--subsets) actions of $\mathbf{G}$ on $\mathbf{U}$ and $\mathbf{%
\check{G}}$ on $\mathbf{\check{U}.}$ We consider an equivariant $\mathbf{%
\check{G}}$ bundle $\mathbf{v}_{\mathbf{U}}\mathbf{\doteqdot U}\times
\mathbf{v,}$ where $\mathbf{v}$ is a d--vector space and a
finite--dimensional representation $\rho :\mathbf{\check{G}\rightarrow }GL(%
\mathbf{v})$ with the $\mathbf{\check{G}}$ action
\begin{equation*}
\mathbf{\check{g}(u,v)=}\left( \mathbf{u\check{g}}^{-1},\rho (\mathbf{\check{%
g}})\mathbf{v}\right) =\left\{
\begin{array}{c}
\left( x\check{g}^{-1},\rho (\check{g})v\right) , \\
\left( y\ ^{\bullet }\check{g}^{-1},\ ^{\bullet }\rho (\ ^{\bullet }\check{g}%
)\ ^{\bullet }v\right)%
\end{array}%
\right\} .
\end{equation*}%
The pull--back of the $\mathbb{Z}_{k}$ distinguished bundle $\mathbf{\check{G%
}\rightarrow G}$ is
\begin{equation*}
\mathbf{\Psi }=\tau ^{\ast }\mathbf{\check{G}\rightarrow U\times U}
\end{equation*}%
defined by $\mathbf{\Psi }_{(\mathbf{u}_{1},\mathbf{u}_{2})}\doteqdot \{\pi (%
\mathbf{\check{g}})=\tau (\mathbf{u}_{1},\mathbf{u}_{2})\}$ for any $\mathbf{%
\check{g}\in \check{G},}$ where $\tau :\mathbf{U\times U\rightarrow }$ $%
\mathbf{G}$ is a canonical map $\mathbf{u}_{1}\tau (\mathbf{u}_{1},\mathbf{u}%
_{2})\rightarrow \mathbf{u}_{2}$ translating $\mathbf{u}_{1}$ into $\mathbf{u%
}_{2}.$ So, $\mathbf{\Psi }_{(\mathbf{u}_{1},\mathbf{u}_{2})}$ is the set of
all distinguished lifts of $\tau (\mathbf{u}_{1},\mathbf{u}_{2})$ to $%
\mathbf{\check{G}}$ and $\mathbf{\Psi }$ is the $\mathbb{Z}_{k}$--principal
bundle provided with a trivial N--connection (in this case, with zero
N--connection curvature). The bundle $\mathbf{\Psi }$ has a module the $%
\mathbf{\check{G}}$--equivariant bundle $\mathbf{v}_{\mathbf{U}}\rightarrow
\mathbf{U.}$ This follows from the fact that for any two pull--backs $%
\mathbf{U}_{1}$ and $\mathbf{U}_{2}$ of $\mathbf{v}_{\mathbf{U}}$ as two
respective projections $\mathbf{U\times U\rightarrow U}$ one has that $%
\mathbf{\Psi }_{(\mathbf{u}_{1},\mathbf{u}_{2})}\subset \mathbf{\check{G}}$
transforms the distinguished fiber in $(\mathbf{u}_{1},\mathbf{u}_{2})$ of $%
\mathbf{U}_{1}$ into the corresponding one of $\mathbf{U}_{2}$ related by
the representation map $\rho .$ Having also the $\mathbf{\check{G}}$%
--equivariance, of $End(\mathbf{v}_{\mathbf{U}}),$ we can write $End(\mathbf{%
v}_{\mathbf{U}})/\mathbf{G=End(v)}$ for distinguished endomorphysms.

\subsubsection{Global constructions}

\label{sgc}The above presented constructions can be globalized to the case
of N--anholonomic manifold $\mathbf{V}$ instead of the d--vector space $%
\mathbf{v}$ (the d--objects with trivial splitting can be considered for any
point of $\mathbf{V}).$ The procedure is completely similar to that given
for ''holonomic'' manifolds in \cite{mursin} but it should be performed in a
form to preserve the N--connection splitting (\ref{whitney}). This may be
achieved by applied globalizing the bundle $\mathbf{\Psi }$ and transforming
it into a nonholonomic bundle.

Having in mind the distinguished extension (\ref{dextseq}), we replace the
set $\mathbf{U}$ by a principal N--anholonomic $\mathbf{G}$--bundle $\pi :%
\mathbf{B\rightarrow V}$ and consider the product $\mathbf{B\times
B\rightarrow V}$ instead of $\mathbf{U\times U.}$ Like for a trivial point
of $\mathbf{U,}$ the globalized map $\tau :\mathbf{B\times B\rightarrow G}$
allows us to introduce $\mathbf{\Psi =}\tau ^{\ast }\mathbf{\check{G}}$
being the $\mathbb{Z}_{k}$--bundle over $\mathbf{B\times B}$ which defines a
lifting N--anholonomic bundle gerbe if to follow the terminology for
holonomic constructions, \cite{murray}.

We can consider d--tensor objects of weight $q$ as $\mathbf{\Psi }^{q}$%
--modules being nonholonomic variants of bundle gerbe modules for the
N--anholonomic bundle gerbe $\mathbf{\Psi }^{q}\doteqdot \mathbf{\Psi }%
^{\otimes q}.$ In more explicit form, we use a $\mathbf{\check{G}}$%
--equivariant bundle $\mathbf{W\rightarrow B}$ for the action of $\mathbf{%
\check{G}}$ of $\mathbf{B:}$

\begin{definition}
The N--anholonomic $\mathbf{\check{G}}$--equivariant bundle $\mathbf{%
W\rightarrow B}$ with defined action of weight $q$ of the isotropy
distinguished subgroups states $\mathbf{W}$ as a $\mathbf{\Psi }^{q}$%
--module.
\end{definition}

The space $\mathbf{W}$ can be also treated as a vector bundle direct sum of $%
\mathbf{\Psi }^{q}$--modules, all adapted to the N--connection structure,
i.e. preserving the h-- and v--decomposition by (\ref{whitney}). This allows
us to concentrate the attention only to ''boldfaced'' $\mathbf{\Psi }^{q}$%
--modules carrying out all information about nonholnomic and non--trivial
topological configurations. Such constructions run parallel to the usual
theory of vector bundles provided with N--connection structure and in a more
formalized form (unifying the approaches to gauge fields, gravity and
geometrized mechanics) to N--anholonomic manifolds.

It should be noted that if $\mathbf{W}$ is a $\mathbf{\Psi }$--module, then
we get a trivial module but it can provided with a nontrivial N--connection
(with nonvanishing N--connection curvature). In such cases, one works with
constructions of type $\mathcal{E}(\mathbf{W})=End(\mathbf{W})/\mathbf{G}$
splitting into h-- and v--subspaces and this allow us to reformulate in
nonholnomic form, for N--anholonomic $\mathbf{\Psi }$--modules, the main
properties of such spaces formally formulated for trivial N--connection
structure \cite{bcmms}.

\begin{proposition}
The $\mathbf{\Psi }^{q}$--modules satisfy the following N--adapted
properties:

\begin{enumerate}
\item N--anholonomic $\mathbf{\Psi }^{q}$--modules and bundles on $\mathbf{V}
$ are bijective equivalent.

\item The bundle of N--adapted endomorphisms of a $\mathbf{\Psi }^{q}$%
--module is a $\mathbf{\Psi }^{0}$--module.

\item The direct sum of two $\mathbf{\Psi }^{q}$--modules is a $\mathbf{\Psi
}^{q}$--module.

\item The d--tensor product of a $\mathbf{\Psi }^{q_{1}}$--module to a $%
\mathbf{\Psi }^{q_{2}}$--module results in a $\mathbf{\Psi }^{q_{1}+q_{2}}$%
--module.
\end{enumerate}
\end{proposition}

We omit the proof of these properties following from an explicit Cech
description of the above structures in N--adapted from (dubbing the
constructions from \cite{bcmms} for h-- and v--configurations). Here we note
that the elements of cohomological classes, like $[e]\in H^{3}(\mathbf{V},%
\mathbb{Z}_{k})\simeq H^{2}(\mathbf{V},U(1))$ and $\delta \lbrack e]\in
H^{3}(\mathbf{V},\mathbb{Z}_{k}),$ are defined for N--anholonomic manifolds,
see Ref. \cite{karoubi} for an introduction in $K$--theory and related
cohomological calculus. This results in distinguished (by N--connection
structure) $K$--group of the semi--group of N--anholonomic $\mathbf{\Psi }$%
--modules. \footnote{%
We emphasize that we wrote $\mathbf{\Psi }$--modules instead of $\mathbf{%
\Gamma }$--modules \cite{murray,bcnns,mursin,aris,bcmms} because in this
work the symbol $\mathbf{\Gamma }$ is used for d--connections.}

\subsection{Curvatures for N--anholonomic bundle ge\-rbe modules}

For a holonomic manifold, because $\mathbb{Z}_{k}$ is finite, there is a
natural $\mathbf{\check{G}}$--equivariant flat connection $\ ^{flat}\nabla _{%
\mathbf{X}}$ on any cart from a covering of bundle $\mathbf{v}_{\mathbf{U}}.$
For N--anholonomic manifolds the role of flat connection is played by metric
canonical d$_{N}$--connection $\widehat{\mathbf{D}}^{N}$ defined by\textbf{\
}a (pseudo) Euclidean N--metric structure \textbf{\ }$^{\eta }\mathbf{g}%
=\eta \oplus _{N}\ ^{\bullet }\eta $ and the N--connection $\mathbf{N,}$ see
Corollary \ref{cncs}. If a d--metric structure $\mathbf{g}=[g,\mathbf{\ }%
^{\bullet }g]$ is stated on such a N--anholonomic manifold,\textbf{\ }we
shall work with the corresponding canonical d--connection $\widehat{\mathbf{D%
}},$ see Theorem \ref{thcdc}. For simplicity, in this section we shall
derive our constructions starting from $\widehat{\mathbf{D}}^{N}$ but we
note that, in general, we can work with an arbitrary d--connection $\mathbf{D%
}$ lifted on $\mathbf{W}$ as a distinguished linear operator, N--adapted to (%
\ref{whitney}), acting in the space of d--forms $\mathbf{\omega =(\omega }%
^{0},\mathbf{\omega }^{1},...),$ where, for instance, $\mathbf{\omega }^{1}$
denotes the space of 1--forms distinguished by the N--connection structure.
Let us consider the d--operator%
\begin{equation*}
\overleftarrow{\mathbf{D}}:\mathbf{\omega }^{0}(\mathbf{B,W})\rightarrow
\mathbf{\omega }^{1}(\mathbf{B,W}).
\end{equation*}%
For a necessary small open subset $\mathbf{U}\subset $ $\mathbf{V},$ we can
identify the restriction of $\mathbf{B}$ to $\mathbf{U}$ with $\mathbf{%
U\times V}$ and their restriction of $\mathbf{W}$ with $\mathbf{W}_{\mathbf{V%
}}.$ In result, we may write
\begin{equation}
\overleftarrow{\mathbf{D}}=\widehat{\mathbf{D}}_{\mathbf{V}}^{N}+\mathbf{D}%
_{B},  \label{dcpsdc}
\end{equation}%
for $\mathbf{D}_{B}$ being a pull--back of a connection from the base $%
\mathbf{U}.$ This way, $\overleftarrow{\mathbf{D}}$ is defined as a $\mathbf{%
\Psi }$--module d--connection if it is equivariant for the group $\mathbf{%
\check{G}}^{\#}\mathbf{=}\left( U(1)\times \mathbf{\check{G}}\right) /%
\mathbb{Z}_{k}$ with $\mathbb{Z}_{k}\subset U(1)\times \mathbf{\check{G}}$
parametrized as a h-- and v--distinguished inclusions by anti--diagonal
subroups. Such a d--connection satisfies the rule
\begin{equation}
\overleftarrow{\mathbf{D}}\left( f\varphi \right) =\mathbf{e}_{\mu }(f)%
\mathbf{e}^{\mu }\otimes _{N}\varphi +f\otimes _{N}\overleftarrow{\mathbf{D}}%
\varphi  \label{aux01}
\end{equation}%
for any function $f$ on $\mathbf{V}$ and section $\varphi $ of $\mathbf{W}$
where $\mathbf{e}_{\mu }$ and $\mathbf{e}^{\mu }$ are N--elongated operators
(\ref{dder}) and (\ref{ddif}).

Let us consider two $\mathbf{\Psi }^{q}$--module d--connections $%
\overleftarrow{\mathbf{D}}_{1}$ and $\overleftarrow{\mathbf{D}}_{2}$ on $%
\mathbf{W}.$ The distorsion $\overleftarrow{\mathbf{P}}=$ $\overleftarrow{%
\mathbf{D}}_{1}-$ $\overleftarrow{\mathbf{D}}_{2}$ is $\mathbf{\check{G}}$%
--equivariant and belongs to the d--vector space $\mathbf{\omega }^{1}(%
\mathbf{B,End(W})),$ this follows from (\ref{aux01}). For any vertical to $%
\mathbf{V}$ d--vector $\mathbf{\lambda ,}$ one holds $\mathbf{\lambda
\rfloor }\overleftarrow{\mathbf{P}}=0.$ This allows us to ''divide'' on $%
\mathbf{G}$ and transform $\overleftarrow{\mathbf{P}}$ into an element of $%
\mathbf{\omega }^{1}(\mathbf{B,}\mathcal{E}(\mathbf{W})),$ i.e. by such
N--adapted distorsions we are able to generate all $\mathbf{\Psi }^{q}$%
--module d--connections starting from (\ref{dcpsdc}). In result, we proved

\begin{proposition}
The set of $\mathbf{\Psi }^{q}$--module d--connections on $\mathbf{W}$ is a
N--distin\-gu\-ished affine space generated by N--adapted distorsions as
elements of $\mathbf{\omega }^{1}(\mathbf{B,}$ $\mathcal{E}(\mathbf{W})).$
\end{proposition}

The curvature of a d--connection \ $\overleftarrow{\mathbf{D}}$ (\ref{dcpsdc}%
) is to be constructed by globalizing the results of Theorem \ref{thr}
(which is a very similar to the proof of the previous Proposition):

\begin{theorem}
The curvature of a $\mathbf{\Psi }^{q}$--module d--connections on $\mathbf{W}
$ descends to define an element $\overleftarrow{\mathbf{R}}\in \mathbf{%
\omega }^{2}(\mathbf{B,}\mathcal{E}(\mathbf{W})).$
\end{theorem}

The d--connection \ (\ref{dcpsdc}) is defined by a N--adapted tensor
product. This extends to a straightforward proof of a corresponding result
for curvature:

\begin{corollary}
For any $\mathbf{\Psi }^{q}$--module d--connections $\overleftarrow{\mathbf{D%
}}$ and $\overleftarrow{\mathbf{D}}^{\prime },$ respectively, on
N--anholonomic $\mathbf{\check{G}}$--equivariant bundles $\mathbf{W}$ and $%
\mathbf{W}^{\prime },$ we can compute the curvature of the d--tensor product
connection $\overleftarrow{\mathbf{D}}$ on $\mathbf{W}$ $\otimes \mathbf{W}%
^{\prime },$%
\begin{equation*}
\mathbf{R}_{B}=\overleftarrow{\mathbf{R}}\otimes 1+1\otimes \overleftarrow{%
\mathbf{R}}^{\prime }\in \mathbf{\omega }^{2}(\mathbf{V,}\mathcal{E}(\mathbf{%
W}\otimes \mathbf{W}^{\prime }))=\mathbf{\omega }^{2}(\mathbf{V,}\mathcal{E}(%
\mathbf{W)}\otimes \mathcal{E}(\mathbf{W}^{\prime })).
\end{equation*}
\end{corollary}

In order to define the (twisted)\ Chern character it is enough to have the
data for a $\mathbf{\Psi }$--module d--connection $\overleftarrow{\mathbf{D}}
$ and its descendent curvature $\overleftarrow{\mathbf{R}}$ \cite%
{bcmms,mursin}. For N--anholonomic configurations, the constructions depend
on the fact if there N--anholonomic manifold is provided or not with a
d--metric structure.

\begin{definition}
The (twisted and nonholonomic) Chern character of a $\mathbf{\Psi }^{q}$%
--mo\-du\-le is defined by the curvature of d--connection $\overleftarrow{%
\mathbf{D}}$ induced by the N--connec\-ti\-on structure,%
\begin{equation}
ch(\overleftarrow{\mathbf{R}})=tr\exp \frac{\overleftarrow{\mathbf{R}}}{2\pi
i}.  \label{chch1}
\end{equation}
\end{definition}

\begin{remark}
If additionally to the N--connection structure on $\mathbf{V,}$ it is
defined a d--metric structure $\mathbf{g,}$ the corresponding Chern
character must be computed by using the $\overleftarrow{\mathbf{R}}^{\prime
} $ defined as a distorsion from the nonholonomic configuration stated by a
d--metric \textbf{\ }$^{\eta }\mathbf{g}=[\eta ,\mathbf{\ }^{\bullet }\eta ]$
(inducing together with $N$ the canonical d--connection $\widehat{\mathbf{D}}%
^{N}$ and $\overleftarrow{\mathbf{R}})\ $to a d--connection $\mathbf{g}=[g,%
\mathbf{\ }^{\bullet }g]$ (inducing the canonical d--connection $\widehat{%
\mathbf{D}}$ and curvature $\overleftarrow{\mathbf{R}}^{\prime }.$
\end{remark}

The values $ch(\overleftarrow{\mathbf{R}})$ and/or $ch(\overleftarrow{%
\mathbf{R}}^{\prime })$ are closed and this mean that the corresponding de
Rham cohomology classes are independent of the choice of $\mathbf{\Psi }^{q}$%
--module d--connections if a N--connection structure is prescribed. This has
a number of interesting applications in modern geometric mechanics,
generalized Finsler geometry and gravity with nontrivial N--anholonomic
structures:

\begin{result}
\label{result3}Any regular Lagrange, or Finsler, configuration is
topologically characterized by the corresponding canonical (twisted) Chern
character (\ref{chch1}) computed by using the curvature $\ ^{L}%
\overleftarrow{\mathbf{R}},$ or $\ ^{F}\overleftarrow{\mathbf{R}},$ induced
by the curvature (\ref{curv1}) defined by the $\ $N--connection$\ ^{L}%
\mathbf{N,}$ or$\ ^{F}\mathbf{N,}$ in the form (\ref{cncl}) and $\ $%
d--metric $\ ^{L}\mathbf{g,}$ or$\ ^{F}\mathbf{g,}$ in the form (\ref{slm})
\ defining the canonical d--connection $\ ^{L}\widehat{\mathbf{D}}$ ($\ ^{F}%
\widehat{\mathbf{D}}).$
\end{result}

The set of exact solutions with generic off--diagonal metrics, nonholnomic
frames and various type of local anisotropy, noncommutative and/or Lie
algebroid symmetries constructed in Refs. \cite%
{vhep2,vs,esv,vncl,vesnc,vesalg,vclalg,vggmaf,vstav} can be globalized for
gerbe configurations with nontrivial N--connection structure, i.e. one holds

\begin{result}
\label{result4}The geometric objects for a N--anholonomic Riemann--Cartan
manifold $\ ^{RC}\mathbf{V}$ can be globalized to N--anholonomic gerbe
configurations and characterized by the corresponding (twisted--anholonomic)
Chern character (\ref{chch1}). This character is computed by using the
curvature $\ \overleftarrow{\mathbf{R}}$ induced by the curvature (\ref%
{curv1}) defined by the $\ $N--connection$\ \mathbf{N,}$ d--metric $\
\mathbf{g}$ (\ref{dmetra}) and the canonical d--connection $\ \widehat{%
\mathbf{D}}.$
\end{result}

Finally, in this section, we conclude that the last two Results state new
types of (topological) symmetries and a new classification of regular
Lagrange systems, Finsler spaces and Einstein--Cartan spaces provided with
N--connection structure.

\section{Nonholonomic Clifford Gerbes and Modules}

This section presents a development of the geometry of N--anholonomic
manifolds and related nonholonomic Clifford and Dirac structures \cite%
{vfs,vhs,vncl}. The reader may consult Refs. \cite{vstav,vv,vclalg} for
local component representations of the results and related local calculus
and proofs.

\subsection{Clifford d--algebras and N--anholonomic bundles}

This work states an explicit example of generalized spinor constructions by
considering in sequence (\ref{dextseq}) the d--groups $\mathbf{G=Spin(}n+m%
\mathbf{)}$ and $\mathbf{\check{G}=SO}$ $(n+m)$ where the boldfaced
d--groups split respectively into h-- and v--components $\mathbf{Spin(}n+m%
\mathbf{)=\{}Spin(n),Spin(m)\mathbf{\}}$ and $\mathbf{SO(}n+m)=\{SO(n),$ $%
SO(m)\}.$ One get the central extension
\begin{equation*}
0\rightarrow \mathbb{Z}_{k}\rightarrow \mathbf{Spin(}n+m\mathbf{)\rightarrow
\mathbf{SO(}n+m\mathbf{)}\rightarrow }1
\end{equation*}%
splitting into respective h-- and v--components,%
\begin{equation*}
0\rightarrow \mathbb{Z}_{k}\rightarrow Spin(n)\mathbf{\rightarrow }SO(n)%
\mathbf{\rightarrow }1
\end{equation*}%
and%
\begin{equation*}
0\rightarrow \mathbb{Z}_{k}\rightarrow Spin(m)\mathbf{\rightarrow }\ SO(m)%
\mathbf{\rightarrow }1.
\end{equation*}

Let us consider two real vector spaces $v$ and $\ ^{\bullet }v$ of dimension
$n$ and $m$ each provided with positive defined scalar products and defining
a d--vector space $\mathbf{v}=v\oplus \ ^{\bullet }v.$ We denote by $C(v)$
and $C(\ ^{\bullet }v)$ the corresponding $\mathbb{Z}_{2}$ graded Clifford
algebras defining a Clifford d--algebra
\begin{equation*}
\mathbf{C}(\mathbf{v})=C_{+}(v)\oplus C_{-}(v)\oplus _{N}C_{+}(\ ^{\bullet
}v)\oplus C_{-}(\ ^{\bullet }v).
\end{equation*}%
The splitting $\pm $ \ is related to the chirality operator $\gamma =\pm $
on $C_{\pm }.$ A hermitian Clifford d--module is a $\mathbb{Z}_{2}$--graded
d--vector space $\mathbf{v}^{E}$ provided with complex scalar products on
the h-- and v--components. The endomorphisms of spin representation $%
S=S^{+}\oplus S^{-}$ and $\ ^{\bullet }S=\ ^{\bullet }S^{+}\oplus \
^{\bullet }S^{-}$ define respectively the hermitian Clifford modules for
conventional h-- and v--subspaces, $C(v)=End(S)$ and $C(\ ^{\bullet
}v)=End(\ ^{\bullet }S).$ Any hermitian Clifford d--modules of finite
dimension can be represented in the form $\mathbf{v}^{E}=S\otimes \mathbf{v}%
^{C}$ where $\mathbf{v}^{C}$ is a complex d--vector space on which $\mathbf{C%
}(\mathbf{v})$ acts trivially in distinguished from. We can identify
\begin{equation*}
\mathbf{v}^{C}=Hom(\mathbf{S,v}^{E})\mbox{ and }End(\mathbf{v}^{C})=End(%
\mathbf{v}^{E})
\end{equation*}%
supposing that such maps are h- and v--split and commute with the action of $%
\mathbf{C}(\mathbf{v}).$

We take the bundle $\mathbf{B}$ from section \ref{sgc} to be the bundle of
N--adapted orthogonal frames (see Proposition \ref{pddv}) on $T\mathbf{V}.$
In the spin case, the construction of $\mathbf{\check{G}}^{\#}$ is that for
the d--group $\mathbf{Spin}_{c}(n+m)$ considered in \cite{vfs,vhs,vncl,vv}
for spin N--anholonomic manifolds, when $\mathbf{\Psi }$ is a trivial $%
\mathbb{Z}_{2}$ N--anholonomic bundle gerbe and $spin$--c when $\mathbf{\Psi
}_{c}$ is a trivial $U(1)$ bundle gerbe provided with N--connection stucture.

\begin{definition}
\label{dsbgb}The lifting bundle gerbe $\mathbf{\Psi }$ for the case $\mathbf{%
G=Spin(}n+m\mathbf{)}$ and $\mathbf{\check{G}=SO(}n+m\mathbf{)}$ is called
the spin--bundle N--anholonomic gerbe.
\end{definition}

We can consider half--spin representations $S^{\pm }$ of $Spin(n)$ and $\
^{\bullet }S^{\pm }$ of $Spin$ $(m)$ and introduce the d--spin
representations
\begin{equation}
\mathbf{S}=\left( S^{+}\oplus S^{-}\right) \oplus _{N}\left( \ ^{\bullet
}S^{+}\oplus \ ^{\bullet }S^{-}\right) .  \label{aux03}
\end{equation}

\begin{definition}
The $\mathbf{\Psi }^{1}$--modules associated to the N--adapted d--spin
representation (\ref{aux03}) \ define the N--anholonomic spin $\mathbf{\Psi }%
^{1}$--modules generalizing the concept of d--spin bundles on $\mathbf{V}.$
\end{definition}

The above mentioned spin constructions have a straightforward extension to
even--dimensional oriented N--anholomic Riemann--Cartan manifolds (this
holds always for oriented Lagrange--Finsler spaces), denoted $\mathbf{V}%
^{2n}.$ One introduces the N--anholonomic bundle of complex Clifford
d--algebras \cite{vncl} of $T^{\ast }\mathbf{V}^{2n}$ and consider the
Clifford distinguished map (multiplication) $\mathbf{c}:T^{\ast }\mathbf{V}%
^{2n}$ $\rightarrow \mathbf{C}(\mathbf{V}^{2n}),$ where formally $\mathbf{%
v\rightarrow V}^{2n}.$

\begin{definition}
A N--anholonomic Clifford module (in brief, Clifford d--mo\-du\-le) is a
complex $\mathbb{Z}_{2}$--graded hermitian N--anholonomic vector bundle
\begin{equation*}
\mathbf{E}=E_{+}\oplus E_{-}\oplus _{N}\ ^{\bullet }E_{+}\oplus \ ^{\bullet
}E_{-}
\end{equation*}%
over $\mathbf{V}^{2n}$ satisfying the properties that $\mathbf{E}_{u}$ is a
hermitian Clifford d--module for $\mathbf{C}_{u}(\mathbf{V}^{2n})$ in each
point $\mathbf{u}\in \mathbf{V}^{2n}$ $\ $and that the sub--bundles $E_{+}$
and $^{\bullet }E_{+}$ are respectively orthogonal to $E_{-}$ and $^{\bullet
}E_{-}.$
\end{definition}

We consider the spin--bundle N--anholonomic gerbe $\mathbf{\Psi }$ from
Definition \ref{dsbgb} and the pull--back of $\mathbf{E}$ to $\mathbf{B},$
denoted $\mathbf{E}_{\mathbf{B}}=\pi ^{-1}(\mathbf{E}),$ where $\pi :\mathbf{%
B}\rightarrow \mathbf{V}^{2n}$ is the bundle of N--adapted frames on $%
\mathbf{V}^{2n}.$ In any point $b\in \mathbf{B,}$ there is an isomorphism
transforming $\left( \mathbf{E}_{\mathbf{B}}\right) _{p}$ into $\mathbf{C}(%
\mathbb{R}^{n+m})$ Clifford d--module. We have $\left( \mathbf{E}_{\mathbf{B}%
}\right) _{p}=\mathbf{S\otimes v}_{(b)}^{C}$ for $\mathbf{v}%
_{(b)}^{C}=Hom\left( \mathbf{S,}\left( \mathbf{E}_{\mathbf{B}}\right)
_{p}\right) $ with the homomorphisms defined on $\mathbf{C}(\mathbb{R}%
^{n+m}).$ The construction can be globalized, $\mathbf{E}_{\mathbf{B}}=%
\mathbf{S\otimes v}^{C}.$ The action of $\mathbf{Spin}(n+m)$ on $\mathbf{S}$
induces also an action $\mathbf{E}_{\mathbf{B}}$ $\ $and transforms it into
a N--anholonomic $\mathbf{\Psi }^{-1}$--module. In result, we proved

\begin{theorem}
For a N--anholonmic spin bundle gerbe $\mathbf{S}_{\mathbf{B}},$ every
Clifford d--module $\mathbf{E}$ on N--anholonomic $\mathbf{V}^{2n},$ with
its nonholonomic bundle gerbe $\mathbf{\Psi }$ , has the form $\mathbf{E}_{%
\mathbf{B}}=\mathbf{S\otimes v}^{C}$ for some N--anholonomic $\mathbf{\Psi }%
^{-1}$--module $\mathbf{v}^{C}.$
\end{theorem}

This theorem generalizes for N--anholonomic spaces some similar results
given in \cite{bgv,mursin}. For spin N--anholonomic manifolds $\mathbf{V}%
^{2n}$ considered in \cite{vfs,vhs,vncl,vv}, we have that every Clifford
d--module is a d--tensor product of an N--anholonomic spin bundle with an
arbitrary bundle.

\subsection{N--anholonomic Dirac operators and gerbes}

\subsubsection{The index topological formula for holonomic Dirac operators}

Let us remember the Atiyah--Singer index topological formula for the Dirac
operator \cite{as1,as2}:%
\begin{equation}
ind(D_{\Gamma }^{+})\doteqdot \dim \ker (D_{\Gamma }^{+})-\dim \ co\ker
(D_{\Gamma }^{+})=<\widehat{A}(M)ch(W),[M]>  \label{indf1}
\end{equation}%
where the compact $M$ is an oriented even dimensional spin manifold with
spin--bundles $S^{\pm }$ and $\Gamma $ is a unitary connection on the vector
bundle $W,$ see details on definitions and denotations in Ref. \cite{bgv}
(below, we shall give details for N--anholonomic configurations). In this
formula, we use the genus of the manifold
\begin{equation*}
\widehat{A}(M)\doteqdot \left| \det \frac{R}{2\sinh (R/2)}\right| ^{1/2}
\end{equation*}%
determined by the Riemannian curvature $R$ of the manifold $M.$ The operator
$D_{\Gamma }^{+}$ is the so--called coupled Dirac operator (first order
differential operator) acting in the form $D_{\Gamma }^{+}:C^{\infty
}(M,E^{+})\rightarrow C^{\infty }(M,E^{-}),$ for $E^{\pm }\doteqdot S^{\pm
}\otimes W.$ This Dirac operator can be introduced for non--spin manifolds
even itself this object is not well defined. In a formal way, we can induce
the Dirac operator as a compatible connection on $E=E^{+}\oplus E^{-}$
treated as a Clifford module with multiplication extended to act as the
identity on $W.$

For non--spin manifolds, one exists an index formula for Dirac operators
defined on hermitian Clifford modules $E,$%
\begin{equation}
ind(D_{\Gamma }^{+})=<\widehat{A}(M)ch(E/S),[M]>  \label{indf2}
\end{equation}%
which, if $M$ is spin and $E^{\pm }\doteqdot S^{\pm }\otimes W,$ the
relative Chern character $ch(E/S)$ reduces to the Chern character of $W,$
i.e. to $ch(W).$ We note that one may be not possible to define a canonical
trivialization of $M$ $\ $but it is supposed that one exist a canonical no
were vanishing (volume) density $[M]$ which allows us to perform the
integration. This always holds for the Riemannian manifolds. The aim of next
section is to prove that formulas (\ref{indf1}) and (\ref{indf2}) can be
correspondingly generalized for N--anholonomic manifolds provided with
d--metric and d--connection structures.

\subsubsection{Twisted nonnholonomic Dirac operators on Clifford gerb\-es}

\label{ssa1}Let us go to the Definition \ref{dsbgb} of the spin--bundle
N--anholonomic gerbe $\mathbf{\Psi }$ derived for an N--anholonmic manifold $%
\mathbf{V.}$ The Clifford multiplication is parametrized by N--adapted maps
between such $\mathbf{\Psi }$--modules,
\begin{equation*}
c:\left( \mathbb{R}_{\mathbf{B}}^{n}\right) ^{\ast }\otimes S_{\mathbf{B}%
}^{+}\rightarrow S_{\mathbf{B}}^{-}\mbox{ and }\ ^{\bullet }c:\left( \mathbb{%
R}_{\mathbf{B}}^{m}\right) ^{\ast }\otimes \ ^{\bullet }S_{\mathbf{B}%
}^{+}\rightarrow \ ^{\bullet }S_{\mathbf{B}}^{-}
\end{equation*}%
where $\mathbb{R}_{\mathbf{B}}^{n+m}=\pi ^{\ast }T\mathbf{V}$ is the
bull--back to the N--adapted frame bundle from the tangent bundle $T\mathbf{V%
}$ with $\mathbb{R}^{n}$ and $\mathbb{R}^{m}$ being the fundamental
representations, respectively, of $SO(n)$ and $SO(m)$ defining the d--group $%
\mathbf{SO}(n+m).$ Any d--connection on $\mathbf{V}$ defines a canonical
d--connection inducing a standard d--connection on the bundle of N--adapted
frames $\mathbf{B}.$

\begin{definition}
\label{twoper}The N--anholonomic (twisted) Dirac operator is defined:%
\begin{equation*}
\mathbb{D}^{+}:C^{\infty }(\mathbf{B},S_{\mathbf{B}}^{+})\rightarrow
C^{\infty }(\mathbf{B},S_{\mathbf{B}}^{-}),
\end{equation*}%
for $\mathbf{\Psi }$--modules and
\begin{equation*}
\overleftarrow{\mathbb{D}}:C^{\infty }(\mathbf{B},S_{\mathbf{B}}^{+}\otimes
\mathbf{W})\rightarrow C^{\infty }(\mathbf{B},S_{\mathbf{B}}^{-}\otimes
\mathbf{W}),
\end{equation*}%
for $\mathbf{W}$ being a $\mathbf{\Psi }^{-1}$--module with induced
canonical d--connection.
\end{definition}

The introduced d--operators split into N--adapted components, $\mathbb{D}%
^{+}=\left( h\mathbb{D}^{+},\ ^{\bullet }\mathbb{D}^{+}\right) $ and $%
\overleftarrow{\mathbb{D}}=\left( h\overleftarrow{\mathbb{D}},\ ^{\bullet }%
\overleftarrow{\mathbb{D}}\right) .$ These operators are correspondingly $%
Spin$ $(n)$-- and $Spin(m)$--invariant. The space $S_{\mathbf{B}}^{\pm
}\otimes \mathbf{W}$ is a N--anholonomic $\mathbf{\Psi }^{0}$--module
descending to bundles $\mathbf{E}^{\pm }$ on\thinspace $\mathbf{V}$ which
transforms the Dirac d--operator to be a twisted Dirac d--operator:%
\begin{equation}
\overleftarrow{\mathbb{D}}^{+}:C^{\infty }(\mathbf{B},\mathbf{E}%
^{+})\rightarrow C^{\infty }(\mathbf{B},\mathbf{E}^{-}).  \label{diropa1}
\end{equation}%
If $\mathbf{V}$ is a spin manifold, than the operators from Definition \ref%
{twoper} descend to $\mathbf{V}$ and with a local decomposition $\mathbf{E}_{%
\mathbf{V}}=\mathbf{S}_{\mathbf{V}}\otimes \mathbf{W}_{\mathbf{V}}.$

Let us consider an N--anholonomic Clifford module $\mathbf{E}$ for $\mathbf{C%
}(\mathbf{V})$ for which there is a d--connection $\ ^{A}\mathbf{D}$ induced
by the canonical d--connection $\widehat{\mathbf{D}}$ and acting following
the rule
\begin{equation*}
\ ^{A}\mathbf{D[c(\omega )}f\mathbf{]=c(\widehat{\mathbf{D}}\omega )\lambda
+c(\omega )}\ ^{A}\mathbf{D\lambda }
\end{equation*}%
for any 1--form $\mathbf{\omega }$ on $\mathbf{V}$ and $\mathbf{\lambda \in }%
C^{\infty }(\mathbf{V,E}).$ This action preserves the global decomposition (%
\ref{whitney}). We may associate to $^{A}\mathbf{D}$ a Dirac d--operator $%
\overleftarrow{\mathbb{D}}$ by using the sequence%
\begin{equation*}
C^{\infty }(\mathbf{V,E})\overset{\ ^{A}\mathbf{D}}{\rightarrow }C^{\infty }(%
\mathbf{V,}T^{\ast }\mathbf{V}\otimes \mathbf{E})\overset{\ \mathbf{c}}{%
\rightarrow }C^{\infty }(\mathbf{V,E}).
\end{equation*}%
This sequence splits also in h-- and v--components. Because the Clifford
multiplication by $T^{\ast }\mathbf{V}$ results in two distinguished odd
parts of $C(\mathbf{V}),$ we get an odd operator $\overleftarrow{\mathbb{D}}$
splitting in $\pm $ components,%
\begin{equation*}
\overleftarrow{\mathbb{D}}^{\pm }:C^{\infty }(\mathbf{V,E}^{\pm
})\rightarrow C^{\infty }(\mathbf{V,E}^{\mp })
\end{equation*}%
acting in distinguished form on h-- and v--components,
\begin{equation*}
\overleftarrow{\mathbb{D}}^{\pm }:C^{\infty }(V\mathbf{,E}^{\pm
})\rightarrow C^{\infty }(V\mathbf{,E}^{\mp })\mbox{ and }\overleftarrow{%
\mathbb{D}}^{\pm }:C^{\infty }(\ ^{\bullet }V\mathbf{,E}^{\pm })\rightarrow
C^{\infty }(\ ^{\bullet }V\mathbf{,E}^{\mp }).
\end{equation*}%
Such operators are N--adapted and formal adjoint of each other respective
h-- and v--component of the standard functional $L^{2}$ (not confusing with
the Lagrange fundamental function considered in the previous section)
defining the inner product on $C^{\infty }(\mathbf{V,E}).$

The curvature $\ ^{A}\mathbf{R}$ of the d--connection $^{A}\mathbf{D}$ is a
2--form with values in
\begin{equation*}
End(\mathbf{E})=\mathbf{C}(\mathbf{V})\otimes End_{\mathbf{C}(\mathbf{V})}(%
\mathbf{E}).
\end{equation*}%
In general, $\ ^{A}\mathbf{R}$ does not commute with the action on $\mathbf{C%
}(\mathbf{V}).$ For Riemannian manifolds, it was proposed to introduce the
twisting curvature \cite{bgv} $R_{E/S}=\ ^{A}R-c(R)$ for any $c(R)\in C(M)$
satisfying the conditions $\left[ \ ^{A}R,c(\lambda )\right] =c(R(\lambda ))$
and $\left[ c(R),c(\lambda )\right] =c(R(\lambda ))$ for $R$ being the
Riemannian curvature of $M$ and any tangent vector $\lambda .$ In a similar
form, \ for N--anholonomic manifolds, we can define the twisting canonical
curvature
\begin{equation*}
\widehat{\mathbf{R}}_{E/S}=\ ^{A}\mathbf{R}-c(\widehat{\mathbf{R}})
\end{equation*}%
induced by\ (\ref{curv1}), see formulas (\ref{dcurv}) from Appendix,
computed for the canonical d--connection $\widehat{\mathbf{D}},$ see Theorem %
\ref{thcdc} and (\ref{candcon}). With this curvature, we may act as with the
Riemannian one following the procedure of defining generalized Dirac
operators from \cite{bgv}.

\subsubsection{Main result and concluding remarks}

The material of previous section \ref{ssa1} consists the proof of

\begin{theorem}
\label{thmr} (Twisted Index formula for N--anholonomic Dirac \newline
operators). If $\mathbf{V}$ is a compact N--anholonomic manifold, then the
Dirac operator $\overleftarrow{\mathbb{D}}^{+}$ (\ref{diropa1}) satisfies
the index formula
\begin{equation*}
ind(\overleftarrow{\mathbb{D}}^{+})=<\widehat{A}(\mathbf{V})ch(\mathbf{W}),[%
\mathbf{V}]>,
\end{equation*}%
where the genus
\begin{equation*}
\widehat{A}(\mathbf{V})\doteqdot \left| \det \frac{\widehat{\mathbf{R}}}{%
2\sinh (\widehat{\mathbf{R}}/2)}\right| ^{1/2}
\end{equation*}%
is determined by the curvature $\widehat{\mathbf{R}}$ of the canonical
d--connection $\widehat{\mathbf{D}}.$
\end{theorem}

This theorem can be stated for certain particular cases of Lagrange, or
Finsler, geometries and their spinor formulation, for instance, with the aim
to locally anisotropic generalization of the so--called $C$--spaces \cite%
{castro,castrop} which will present topological characteristics derived from
a fundamental Lagrange (or Finsler) function or, in a new fashion, for
non--spin $C$--gerbes associated to nonholonomic gravitational and spinor
interactions. The Main Result of this work can be also applied for
topological classification of new types of globalized exact solutions
defining nonholonomic gravitational and matter field configurations \cite%
{vhep2,vs,esv,vncl,vesnc,vesalg,vclalg,vggmaf,vstav} .

\vskip10pt

\textbf{Acknowledgement: } The authors are grateful to C. Castro Perelman
for useful discussions. S. V. thanks Prof. M. Anastasiei for kind support.

\setcounter{equation}{0} \renewcommand{\theequation}
{A.\arabic{equation}} \setcounter{subsection}{0}
\renewcommand{\thesubsection}
{A.\arabic{subsection}}

\appendix

\section{ Some Local Formulas from N--Connection Ge\-ometry}

In this Appendix, we present some component formulas and equations defining
the local geometry of N--anholonomic spaces, see details in Refs. \cite%
{vggmaf,vncl,vesalg,ma}.

Locally, a N--connection, see Definition \ref{dnc}, is stated by its
coefficients $N_{i}^{a}(u),$
\begin{equation}
\mathbf{N}=N_{i}^{a}(u)dx^{i}\otimes \partial _{a}  \label{nclf}
\end{equation}%
where the local coordinates (in general, abstract ones both for holonomic
and nonholonomic variables) are split in the form $u=(x,y),$ or $u^{\alpha
}=\left( x^{i},y^{a}\right) ,$ where $i,j,k,\ldots =1,2,\ldots ,n$ and $%
a,b,c,\ldots =n+1,n+2,\ldots ,n+m$ when $\partial _{i}=\partial /\partial
x^{i}$ and $\partial _{a}=\partial /\partial y^{a}.$ The well known class of
linear connections consists on a particular subclass with the coefficients
being linear on $y^{a},$ i.e., $N_{i}^{a}(u)=\Gamma _{bj}^{a}(x)y^{b}.$

An explicit local calculus allows us to write the N--connection curvature (%
\ref{njht}) in the form
\begin{equation*}
\mathbf{\Omega }=\frac{1}{2}\Omega _{ij}^{a}dx^{i}\wedge dx^{j}\otimes
\partial _{a},
\end{equation*}%
with the N--connection curvature coefficients
\begin{equation}
\Omega _{ij}^{a}=\delta _{\lbrack j}N_{i]}^{a}=\delta _{j}N_{i}^{a}-\delta
_{i}N_{j}^{a}=\partial _{j}N_{i}^{a}-\partial
_{i}N_{j}^{a}+N_{i}^{b}\partial _{b}N_{j}^{a}-N_{j}^{b}\partial
_{b}N_{i}^{a}.  \label{ncurv}
\end{equation}

Any N--connection $\mathbf{N}=N_{i}^{a}(u)$ induces a N--adapted frame
(vielbein) structure
\begin{equation}
\mathbf{e}_{\nu }=\left( e_{i}=\partial _{i}-N_{i}^{a}(u)\partial
_{a},e_{a}=\partial _{a}\right) ,  \label{dder}
\end{equation}%
and the dual frame (coframe) structure%
\begin{equation}
\mathbf{e}^{\mu }=\left( e^{i}=dx^{i},e^{a}=dy^{a}+N_{i}^{a}(u)dx^{i}\right)
.  \label{ddif}
\end{equation}%
The vielbeins (\ref{ddif}) satisfy the nonholonomy (equivalently,
anholonomy) relations
\begin{equation}
\lbrack \mathbf{e}_{\alpha },\mathbf{e}_{\beta }]=\mathbf{e}_{\alpha }%
\mathbf{e}_{\beta }-\mathbf{e}_{\beta }\mathbf{e}_{\alpha }=W_{\alpha \beta
}^{\gamma }\mathbf{e}_{\gamma }  \label{anhrel}
\end{equation}%
with (antisymmetric) nontrivial anholonomy coefficients $W_{ia}^{b}=\partial
_{a}N_{i}^{b}$ and $W_{ji}^{a}=\Omega _{ij}^{a}.$\footnote{%
One preserves a relation to our previous denotations \cite{vfs,vhs} if we
consider that $\mathbf{e}_{\nu }=(e_{i},e_{a})$ and $\mathbf{e}^{\mu
}=(e^{i},e^{a})$ are, respectively, the former $\delta _{\nu }=\delta
/\partial u^{\nu }=(\delta _{i},\partial _{a})$ and $\delta ^{\mu }=\delta
u^{\mu }=(d^{i},\delta ^{a});$ we emphasize that operators (\ref{dder}) and (%
\ref{ddif}) define, correspondingly, the ``N--elongated'' partial
derivatives and differentials which are convenient for calculations on
N--anholonomic manifolds.} These formulas present a local proof of
Proposition \ref{pddv} when
\begin{equation*}
\mathbf{e}=\{\mathbf{e}_{\nu }\}=(\ e=\{e_{i}\},^{\bullet }e=\{e_{a}\})
\end{equation*}%
and
\begin{equation*}
\widetilde{\mathbf{e}}=\{\mathbf{e}^{\mu }\}=\left( \widetilde{e}%
=\{e^{i}\},\ ^{\bullet }\widetilde{e}=\{e^{a}\}\right) .
\end{equation*}

Let us consider metric structure%
\begin{equation}
\ \breve{g}=\underline{g}_{\alpha \beta }\left( u\right) du^{\alpha }\otimes
du^{\beta }  \label{metr}
\end{equation}%
defined with respect to a local coordinate basis $du^{\alpha }=\left(
dx^{i},dy^{a}\right) $ by coefficients%
\begin{equation}
\underline{g}_{\alpha \beta }=\left[
\begin{array}{cc}
g_{ij}+N_{i}^{a}N_{j}^{b}h_{ab} & N_{j}^{e}h_{ae} \\
N_{i}^{e}h_{be} & h_{ab}%
\end{array}%
\right] .  \label{ansatz}
\end{equation}%
In general, such a metric (\ref{ansatz})\ is generic off--diagonal, i.e. it
can not be diagonalized by any coordinate transforms and that $N_{i}^{a}(u)$
are any general functions. The condition (\ref{algn01}), for $X\rightarrow
e_{i}$ and $\ ^{\bullet }Y\rightarrow \ ^{\bullet }e_{a},$ transform into
\begin{equation*}
\breve{g}(e_{i},\ ^{\bullet }e_{a})=0,\mbox{ equivalently }\underline{g}%
_{ia}-N_{i}^{b}h_{ab}=0
\end{equation*}%
where $\underline{g}_{ia}$ $\doteqdot g(\partial /\partial x^{i},\partial
/\partial y^{a}),$ which allows us to define in a unique form the
coefficients $N_{i}^{b}=h^{ab}\underline{g}_{ia}$ where $h^{ab}$ is inverse
to $h_{ab}.$ We can write the metric $\breve{g}$ with ansatz (\ref{ansatz})\
in equivalent form, as a d--metric adapted to a N--connection structure, see
Definition \ref{ddm},
\begin{equation}
\mathbf{g}=\mathbf{g}_{\alpha \beta }\left( u\right) \mathbf{e}^{\alpha
}\otimes \mathbf{e}^{\beta }=g_{ij}\left( u\right) e^{i}\otimes
e^{j}+h_{ab}\left( u\right) \ ^{\bullet }e^{a}\otimes \ ^{\bullet }e^{b},
\label{dmetr}
\end{equation}%
where $g_{ij}\doteqdot \mathbf{g}\left( e_{i},e_{j}\right) $ and $%
h_{ab}\doteqdot \mathbf{g}\left( \ ^{\bullet }e_{a},\ ^{\bullet
}e_{b}\right) $ \ and the vielbeins $\mathbf{e}_{\alpha }$ and $\mathbf{e}%
^{\alpha }$ are respectively of type (\ref{dder}) and (\ref{ddif}).

We can say that the metric $\ \breve{g}$ (\ref{metr}) is equivalently
transformed into (\ref{dmetr}) \ by performing a frame (vielbein) transform
\begin{equation*}
\mathbf{e}_{\alpha }=\mathbf{e}_{\alpha }^{\ \underline{\alpha }}\partial _{%
\underline{\alpha }}\mbox{ and }\mathbf{e}_{\ }^{\beta }=\mathbf{e}_{\
\underline{\beta }}^{\beta }du^{\underline{\beta }}.
\end{equation*}%
with coefficients

\begin{eqnarray}
\mathbf{e}_{\alpha }^{\ \underline{\alpha }}(u) &=&\left[
\begin{array}{cc}
e_{i}^{\ \underline{i}}(u) & N_{i}^{b}(u)e_{b}^{\ \underline{a}}(u) \\
0 & e_{a}^{\ \underline{a}}(u)%
\end{array}%
\right] ,  \label{vt1} \\
\mathbf{e}_{\ \underline{\beta }}^{\beta }(u) &=&\left[
\begin{array}{cc}
e_{\ \underline{i}}^{i\ }(u) & -N_{k}^{b}(u)e_{\ \underline{i}}^{k\ }(u) \\
0 & e_{\ \underline{a}}^{a\ }(u)%
\end{array}%
\right] ,  \label{vt2}
\end{eqnarray}%
being linear on $N_{i}^{a}.$ We can consider that a N--anholonomic manifold $%
\mathbf{V}$ provided with metric structure $\breve{g}$ (\ref{metr})
(equivalently, with d--metric (\ref{dmetr})) is a special type of a manifold
provided with a global splitting into conventional ``horizontal'' and
``vertical'' subspaces (\ref{whitney}) induced by the ``off--diagonal''
terms $N_{i}^{b}(u)$ and a prescribed type of nonholonomic frame structure (%
\ref{anhrel}).

A d--connection, see Definition \ref{ddc}, splits into h-- and v--covariant
derivatives, $\mathbf{D}=D+\ ^{\bullet }D,$ where $D_{k}=\left(
L_{jk}^{i},L_{bk\;}^{a}\right) $ and $\ \ ^{\bullet }D_{c}=\left(
C_{jk}^{i},C_{bc}^{a}\right) $ are correspondingly introduced as h- and
v--parametrizations of (\ref{cond1}),%
\begin{equation*}
L_{jk}^{i}=\left( \mathbf{D}_{k}e_{j}\right) \rfloor e^{i},\quad
L_{bk}^{a}=\left( \mathbf{D}_{k}e_{b}\right) \rfloor
e^{a},~C_{jc}^{i}=\left( \mathbf{D}_{c}e_{j}\right) \rfloor e^{i},\quad
C_{bc}^{a}=\left( \mathbf{D}_{c}e_{b}\right) \rfloor e^{a}.
\end{equation*}%
The components $\mathbf{\Gamma }_{\ \alpha \beta }^{\gamma }=\left(
L_{jk}^{i},L_{bk}^{a},C_{jc}^{i},C_{bc}^{a}\right) $ completely define a
d--connecti\-on $\mathbf{D}$ on a N--anholonomic manifold $\mathbf{V}.$

The simplest way to perform a local covariant calculus by applying
d--connecti\-ons is to use N--adapted differential forms like $\mathbf{%
\Gamma }_{\beta }^{\alpha }=\mathbf{\Gamma }_{\beta \gamma }^{\alpha }%
\mathbf{e}^{\gamma }$ with the coefficients defined with respect to (\ref%
{ddif}) and (\ref{dder}).One can introduce the d--connection 1--form%
\begin{equation*}
\mathbf{\Gamma }_{\ \beta }^{\alpha }=\mathbf{\Gamma }_{\ \beta \gamma
}^{\alpha }\mathbf{e}^{\gamma },
\end{equation*}%
when the N--adapted components of d-connection $\mathbf{D}_{\alpha }=(%
\mathbf{e}_{\alpha }\rfloor \mathbf{D})$ are computed following formulas
\begin{equation}
\mathbf{\Gamma }_{\ \alpha \beta }^{\gamma }\left( u\right) =\left( \mathbf{D%
}_{\alpha }\mathbf{e}_{\beta }\right) \rfloor \mathbf{e}^{\gamma },
\label{cond1}
\end{equation}%
where ''$\rfloor "$ denotes the interior product. This allows us to define
in local form the torsion (\ref{tors1}) $\mathbf{T=\{\mathcal{T}^{\alpha }\},%
}$ where
\begin{equation}
\mathcal{T}^{\alpha }\doteqdot \mathbf{De}^{\alpha }=d\mathbf{e}^{\alpha
}+\Gamma _{\ \beta }^{\alpha }\wedge \mathbf{e}^{\beta }  \label{tors}
\end{equation}%
and curvature (\ref{curv1}) $\mathbf{R}=\{\mathcal{R}_{\ \beta }^{\alpha
}\}, $ where
\begin{equation}
\mathcal{R}_{\ \beta }^{\alpha }\doteqdot \mathbf{D\Gamma }_{\beta }^{\alpha
}=d\mathbf{\Gamma }_{\beta }^{\alpha }-\Gamma _{\ \beta }^{\gamma }\wedge
\mathbf{\Gamma }_{\ \gamma }^{\alpha }.  \label{curv}
\end{equation}

The d--torsions components of a d--connection $\mathbf{D},$ see Theorem \ref%
{tht}, $\ $are computed
\begin{eqnarray}
T_{\ jk}^{i} &=&L_{\ jk}^{i}-L_{\ kj}^{i},\ T_{\ ja}^{i}=-T_{\ aj}^{i}=C_{\
ja}^{i},\ T_{\ ji}^{a}=\Omega _{\ ji}^{a},\   \notag \\
T_{\ bi}^{a} &=&T_{\ ib}^{a}=\frac{\partial N_{i}^{a}}{\partial y^{b}}-L_{\
bi}^{a},\ T_{\ bc}^{a}=C_{\ bc}^{a}-C_{\ cb}^{a}.  \label{dtors}
\end{eqnarray}%
For instance, $T_{\ jk}^{i}$ and $T_{\ bc}^{a}$ are respectively the
coefficients of the $h(hh)$--torsion $T(X,Y)$ and $v(vv)$--torsion $\mathbf{%
\ }^{\bullet }T(\mathbf{\ }^{\bullet }X,\mathbf{\ }^{\bullet }Y).$

The Levi--Civita linear connection $\nabla =\{^{\nabla }\mathbf{\Gamma }%
_{\beta \gamma }^{\alpha }\},$ with vanishing both torsion and nonmetricity $%
\nabla \breve{g}=0,$ is not adapted to the global splitting (\ref{whitney}).
There is a preferred, canonical d--connection structure,$\ \widehat{\mathbf{D%
}}\mathbf{,}$ $\ $on N--anholonomic manifold $\mathbf{V}$ constructed only
from the metric and N--con\-nec\-ti\-on coefficients $%
[g_{ij},h_{ab},N_{i}^{a}]$ and satisfying the conditions $\widehat{\mathbf{D}%
}\mathbf{g}=0$ and $\widehat{T}_{\ jk}^{i}=0$ and $\widehat{T}_{\ bc}^{a}=0,$
see Theorem \ref{thcdc}. By straightforward calculations with respect to the
N--adapted bases (\ref{ddif}) and (\ref{dder}), we can verify that the
connection
\begin{equation}
\widehat{\mathbf{\Gamma }}_{\beta \gamma }^{\alpha }=\ ^{\nabla }\mathbf{%
\Gamma }_{\beta \gamma }^{\alpha }+\ \widehat{\mathbf{P}}_{\beta \gamma
}^{\alpha }  \label{cdc}
\end{equation}%
with the deformation d--tensor \footnote{$\widehat{\mathbf{P}}_{\beta \gamma
}^{\alpha }$ is a tensor field of type (1,2). As is well known, the sum of a
linear connection and a tensor field of type (1,2) is a new linear
connection.}
\begin{equation}
\widehat{\mathbf{P}}_{\beta \gamma }^{\alpha
}=(P_{jk}^{i}=0,P_{bk}^{a}=e_{b}(N_{k}^{a}),P_{jc}^{i}=-\frac{1}{2}%
g^{ik}\Omega _{\ kj}^{a}h_{ca},P_{bc}^{a}=0)  \label{cdca}
\end{equation}%
satisfies the conditions of the mentioned Theorem. It should be noted that,
in general, the components $\widehat{T}_{\ ja}^{i},\ \widehat{T}_{\ ji}^{a}$
and $\widehat{T}_{\ bi}^{a}$ are not zero. This is an anholonomic frame (or,
equivalently, off--diagonal metric) effect. With respect to the N--adapted
frames, the coefficients\newline
$\widehat{\mathbf{\Gamma }}_{\ \alpha \beta }^{\gamma }=\left( \widehat{L}%
_{jk}^{i},\widehat{L}_{bk}^{a},\widehat{C}_{jc}^{i},\widehat{C}%
_{bc}^{a}\right) $ are computed:
\begin{eqnarray}
\widehat{L}_{jk}^{i} &=&\frac{1}{2}g^{ir}\left(
e_{k}g_{jr}+e_{j}g_{kr}-e_{r}g_{jk}\right) ,  \label{candcon} \\
\widehat{L}_{bk}^{a} &=&e_{b}(N_{k}^{a})+\frac{1}{2}h^{ac}\left(
e_{k}h_{bc}-h_{dc}\ e_{b}N_{k}^{d}-h_{db}\ e_{c}N_{k}^{d}\right) ,  \notag \\
\widehat{C}_{jc}^{i} &=&\frac{1}{2}g^{ik}e_{c}g_{jk},\ \widehat{C}_{bc}^{a}=%
\frac{1}{2}h^{ad}\left( e_{c}h_{bd}+e_{c}h_{cd}-e_{d}h_{bc}\right) .  \notag
\end{eqnarray}

In some approaches to Finsler geometry \cite{cbs}, one uses the so--called
Berwald d--connection $\mathbf{D}^{B}$ with the coefficients
\begin{equation}
\ ^{B}\mathbf{\Gamma }_{\ \alpha \beta }^{\gamma }=\left( \ ^{B}L_{jk}^{i}=%
\widehat{L}_{jk}^{i},\ ^{B}L_{bk}^{a}=e_{b}(N_{k}^{a}),\ ^{B}C_{jc}^{i}=0,\
^{B}C_{bc}^{a}=\widehat{C}_{bc}^{a}\right) .  \label{dbc}
\end{equation}%
This d--connection minimally extends the Levi--Civita connection (it is just
the Levi--Civita connection if the integrability conditions are satisfied,
i.e. $\Omega _{\ kj}^{a}=0,$ see (\ref{cdca})). But, in general, for this
d--connection the metricity conditions are not satisfied, for instance $%
D_{a}g_{ij}\neq 0$ and $D_{i}h_{ab}\neq 0.$

By a straightforward d--form calculus in (\ref{curv}), we can find the
N--adapted components $\mathbf{R}_{\ \beta \gamma \delta }^{\alpha }$ of the
curvature $\mathbf{R=\{\mathcal{R}_{\ \beta }^{\alpha }\}}$ of a
d--connection $\mathbf{D},$ i.e. the d--curvatures from Theorem \ref{thr}:
\begin{eqnarray}
R_{\ hjk}^{i} &=&e_{k}L_{\ hj}^{i}-e_{j}L_{\ hk}^{i}+L_{\ hj}^{m}L_{\
mk}^{i}-L_{\ hk}^{m}L_{\ mj}^{i}-C_{\ ha}^{i}\Omega _{\ kj}^{a},  \notag \\
R_{\ bjk}^{a} &=&e_{k}L_{\ bj}^{a}-e_{j}L_{\ bk}^{a}+L_{\ bj}^{c}L_{\
ck}^{a}-L_{\ bk}^{c}L_{\ cj}^{a}-C_{\ bc}^{a}\Omega _{\ kj}^{c},  \notag \\
R_{\ jka}^{i} &=&e_{a}L_{\ jk}^{i}-D_{k}C_{\ ja}^{i}+C_{\ jb}^{i}T_{\
ka}^{b},  \label{dcurv} \\
R_{\ bka}^{c} &=&e_{a}L_{\ bk}^{c}-D_{k}C_{\ ba}^{c}+C_{\ bd}^{c}T_{\
ka}^{c},  \notag \\
R_{\ jbc}^{i} &=&e_{c}C_{\ jb}^{i}-e_{b}C_{\ jc}^{i}+C_{\ jb}^{h}C_{\
hc}^{i}-C_{\ jc}^{h}C_{\ hb}^{i},  \notag \\
R_{\ bcd}^{a} &=&e_{d}C_{\ bc}^{a}-e_{c}C_{\ bd}^{a}+C_{\ bc}^{e}C_{\
ed}^{a}-C_{\ bd}^{e}C_{\ ec}^{a}.  \notag
\end{eqnarray}

Contracting respectively the components of (\ref{dcurv}), one proves

The Ricci tensor $\mathbf{R}_{\alpha \beta }\doteqdot \mathbf{R}_{\ \alpha
\beta \tau }^{\tau }$ is characterized by h- v--components, i.e. d--tensors,%
\begin{equation}
R_{ij}\doteqdot R_{\ ijk}^{k},\ \ R_{ia}\doteqdot -R_{\ ika}^{k},\
R_{ai}\doteqdot R_{\ aib}^{b},\ R_{ab}\doteqdot R_{\ abc}^{c}.
\label{dricci}
\end{equation}%
It should be noted that this tensor is not symmetric for arbitrary
d--connecti\-ons $\mathbf{D}.$

The scalar curvature of a d--connection is
\begin{equation}
\ ^{s}\mathbf{R}\doteqdot \mathbf{g}^{\alpha \beta }\mathbf{R}_{\alpha \beta
}=g^{ij}R_{ij}+h^{ab}R_{ab},  \label{sdccurv}
\end{equation}%
defined by a sum the h-- and v--components of (\ref{dricci}) and d--metric (%
\ref{dmetr}).

The Einstein tensor is defined and computed in standard form
\begin{equation}
\mathbf{G}_{\alpha \beta }=\mathbf{R}_{\alpha \beta }-\frac{1}{2}\mathbf{g}%
_{\alpha \beta }\ ^{s}\mathbf{R}  \label{enstdt}
\end{equation}

For a Lagrange geometry, see Definition \ref{dlg}, by straightforward
component calculations, one can be proved the fundamental results:

\begin{enumerate}
\item The Euler--Lagrange equations%
\begin{equation*}
\frac{d}{d\tau }\left( \frac{\partial L}{\partial y^{i}}\right) -\frac{%
\partial L}{\partial x^{i}}=0
\end{equation*}%
where $y^{i}=\frac{dx^{i}}{d\tau }$ for $x^{i}(\tau )$ depending on
parameter $\tau ,$ are equivalent to the ``nonlinear'' geodesic equations
\begin{equation*}
\frac{d^{2}x^{i}}{d\tau ^{2}}+2G^{i}(x^{k},\frac{dx^{j}}{d\tau })=0
\end{equation*}%
defining paths of the canonical semispray%
\begin{equation*}
S=y^{i}\frac{\partial }{\partial x^{i}}-2G^{i}(x,y)\frac{\partial }{\partial
y^{i}}
\end{equation*}%
where
\begin{equation*}
2G^{i}(x,y)=\frac{1}{2}\ ^{L}g^{ij}\left( \frac{\partial ^{2}L}{\partial
y^{i}\partial x^{k}}y^{k}-\frac{\partial L}{\partial x^{i}}\right)
\end{equation*}%
with $^{L}g^{ij}$ being inverse to (\ref{lqf}).

\item There exists on $\widetilde{TM}$ a canonical N--connection $\ $%
\begin{equation}
\ ^{L}N_{j}^{i}=\frac{\partial G^{i}(x,y)}{\partial y^{i}}  \label{cncl}
\end{equation}%
defined by the fundamental Lagrange function $L(x,y),$ which prescribes
nonholonomic frame structures of type (\ref{dder}) and (\ref{ddif}), $\ ^{L}%
\mathbf{e}_{\nu }=(e_{i},\ ^{\bullet }e_{k})$ and $\ ^{L}\mathbf{e}^{\mu
}=(e^{i},\ ^{\bullet }e^{k}).$

\item The canonical N--connection (\ref{cncl}), defining $\ \ ^{\bullet
}e_{i},$ induces naturally an almost complex structure $\mathbf{F}:\chi (%
\widetilde{TM})\rightarrow \chi (\widetilde{TM}),$ where $\chi (\widetilde{TM%
})$ denotes the module of vector fields on $\widetilde{TM},$%
\begin{equation*}
\mathbf{F}(e_{i})=\ \ ^{\bullet }e_{i}\mbox{ and }\mathbf{F}(\ \ ^{\bullet
}e_{i})=-e_{i},
\end{equation*}%
when
\begin{equation}
\mathbf{F}=\ \ ^{\bullet }e_{i}\otimes e^{i}-e_{i}\otimes \ \ ^{\bullet
}e^{i}  \label{acs1}
\end{equation}%
satisfies the condition $\mathbf{F\rfloor \ F=-I,}$ i. e. $F_{\ \ \beta
}^{\alpha }F_{\ \ \gamma }^{\beta }=-\delta _{\gamma }^{\alpha },$ where $%
\delta _{\gamma }^{\alpha }$ is the Kronecker symbol and ``$\mathbf{\rfloor }
$'' denotes the interior product.

\item On $\widetilde{TM},$ there is a canonical metric structure%
\begin{equation}
\ ^{L}\mathbf{g}=\ ^{L}g_{ij}(x,y)\ e^{i}\otimes e^{j}+\ ^{L}g_{ij}(x,y)\ \
\ ^{\bullet }e^{i}\otimes \ \ ^{\bullet }e^{j}  \label{slm}
\end{equation}%
constructed as a Sasaki type lift from $M.$

\item There is also a canonical d--connection structure $^{L}\widehat{%
\mathbf{\Gamma }}_{\ \alpha \beta }^{\gamma }$ defined only by the
components of $^{L}N_{j}^{i}$ and $^{L}g_{ij},$ i.e. by the \ coefficients
of metric (\ref{slm}) which in its turn is induced by a regular Lagrangian.
\ The values $\ ^{L}\widehat{\mathbf{\Gamma }}_{\ \alpha \beta }^{\gamma
}=(\ ^{L}\widehat{L}_{jk}^{i},\ ^{L}\widehat{C}_{bc}^{a})$ are computed just
as similar values from (\ref{candcon}). We note that on $\widetilde{TM}$
there are couples of distinguished sets of h- and v--components.
\end{enumerate}

\end{document}